\documentclass{article}[12pt]
\usepackage{amssymb}
\usepackage{latexsym}
\usepackage{pb-diagram}
\usepackage[all]{xy}
\usepackage[mathscr]{euscript}
\usepackage{graphicx}

\newcommand\ul{\underline}

\newcommand\bb{\Bbb}
\newcommand\pr{^{\prime}}

\newcommand\commentout[1]{\marginpar{\tiny $\backslash$commentout}}

\newcommand\qed{\hfill$\square$}

\def\column#1#2{\mathrel{\mathop{#1}\limits_{#2}}}

\def\compcirc {\mbox{\hspace{.05cm}}\raisebox{.04cm}{\tiny  {$\circ$ }}}

\newtheorem{Lemma}{Lemma}[section]
\newtheorem{Theorem}{Theorem}[section]
\newtheorem{Proposition}{Proposition}[section]
\newtheorem{Definition}{Definition}[section]
\newtheorem{Corollary}{Corollary}[section]

\newtheorem{Example}{Example}[section]
\newtheorem{Remark}{Remark}[section]

\newenvironment{Proof}{\par\noindent\textbf{Proof:}}
{\qed}

\usepackage{graphics}
\usepackage[tocflat]{tocstyle}

\title{Topological Methods for the Analysis of High Dimensional Data Sets and 3D Object Recognition}

\title{Persistent Homology and Applied Homotopy Theory}
\author{Gunnar Carlsson \\Mathematics Department \\Stanford University \\ and \\Ayasdi, Inc. }
\makeatletter
\renewcommand\tableofcontents{%
  \null\hfill\textbf{\Large\contentsname}\hfill\null\par
  \@mkboth{\MakeUppercase\contentsname}{\MakeUppercase\contentsname}%
  \@starttoc{toc}%
}
\makeatother

\begin{document}

\maketitle
\tableofcontents
\newpage
\section{Introduction} Persistent homology is a technique that has been developed over the last 20 years. Initial ideas developed in the early 1990's \cite{frosini}, but  the idea of persistence was introduced by Vanessa Robins in \cite{robins},  rapidly followed by additional development (\cite{letscher}, \cite{zomorodian}), and  has been developing rapidly since that time.   The original motivation for the method was  to extend the ideas of algebraic topology from the category of spaces $X$ to situations where we only have a sampling of the space $X$.  Of course, a sample is a discrete space so there is nothing to be obtained unless one retains some additional information. One assumes the presence of a metric or  a more relaxed ``dissimilarity measure", and  uses this  information restricted to the sample in constructing the algebraic invariant.  Over time, persistent homology has been used in other situations, for example where one has a topological space with additional information, such as a continuous real valued function, and the sublevel sets of the function determine a filtration on the space.  The output of standard persistent homology (we will discuss some generalizations) is represented in two ways, via {\em persistence barcodes} and {\em persistence diagrams}.  Initially persistent homology  was used, as homology is used for topological spaces, to obtain a large scale geometric understanding of complex data sets, encoded as finite metric spaces.  Examples of this kind of application are \cite{klein}, \cite{giusti}, \cite{hess}, \cite{hess2}, and \cite{evolution}.  Another class of applications uses persistent homology to study data sets where the {\em points themselves} are  metric spaces, such as  databases of molecule structures or images.  This second set of applications is developing very rapidly, and is exemplified in \cite{guowei1}, \cite{guowei2}, \cite{guowei3}, and \cite{hiraoka}.  Another direction in which persistent homology is being applied is in the study of coverage and evasion problems arising in sensor net technology \cite{sensor}.  Example of research in this direction are \cite{desilva1}, \cite{desilva2}, \cite{evasionadams}, and \cite{ghrist}. 

There are a number of different active research directions in this area.  

\begin{itemize}
\item{{\bf Coordinatization of barcodes:} Barcodes in their native form, i.e. as a set of intervals, do not lend themselves to analysis by machine learning techniques.  It is therefore important to represent them in a method more amenable to analysis.  This can be achieved by appropriate coordinatizations of  the space of barcodes. Several methods for this task are described in Section \ref{coordinatization}.}
\item{{\bf Generalized persistence:} Persistent homology has as its output a diagram of complexes, parametrized by the partially ordered set of real numbers, on which algebraic computations are performed so as to produce barcodes.  There are other parameter categories that are useful in the study of data sets.  We discuss two examples of this notion in Section \ref{generalized}, but we would expect there to be many different types of diagrams that will shed light on finite metric spaces.   }
\item{{\bf Stability results:} Since noise and error are key elements in the study of data, it is important to develop methods that quantify the dependence of the barcode output on small perturbations of data.  This requires the imposition of metrics on the set of barcodes, and proving theorems concerning the distances between barcodes that differ only by small perturbations.  The progress that has been made in this direction is described in Section \ref{wholemetrics}.}
\item{{\bf Probabilistic analysis and inference:} Because many of the applications of persistent homology occur in the study of data, it becomes important to study not only stability results under small perturbations, but also perturbations that are {``probabilistically small"}, i.e. which may be large, but where a large perturbation is a rare event.  This means that one must study the distributions on barcode space that occur from applying persistent homology to complexes generated by various random models.  This is a rapidly developing area within the subject.}
\item{{\bf Coverage and evasion problems:}  This work centers around attempts to understand complements in Euclidean space of regions defined by collections of sensors.  It has been approached using different methods, and appears to be a place where techniques such as Spanier-Whitehead duality and embedding analysis, applied in suitably generalized situations, should play a role.  The general problem of understanding complements of objects embedded in Euclidean space of course also plays a role in robotics.  }
\item{{\bf Symplectic geometry:} Although persistent homology is mainly used in situations where one is examining discrete approximations to continuous objects, it can be applied in any geometric situation where there is a metric, or where one is considering filtered objects.  Such situations occur in symplectic geometry, and there is recent work applying the technique in studying, for example, Floer homology spectra.  Examples of this kind of work are \cite{buhovsky}, \cite{leroux}, \cite{polterovich}, \cite{polterovich1}, and \cite{usher}.  }
\end{itemize}

The goal of this paper is to discuss the different research directions and applications at a high level, so that the reader can orient him/herself in the techniques.  We remark that there are a number of useful surveys on persistent homology and on topological data analysis more generally.  The papers \cite{topdata}, \cite{hareredelsbrunner}, and \cite{otter} give different perspectives on this subject.

\section{The Motivating Problem}

Suppose that we are given a set of points $X$  in the plane, and believe that it is reasonable to assume that the points are sampled (perhaps with error) from a geometric object.  We could ask for information about the homology of the underlying space from which $X$ is sampled.  
Consider the  finite set of points $X$ in $\Bbb{R}^2$ displayed in Figure \ref{statisticcircle} below.

{\begin{figure}[!htp]
\centering
\includegraphics[width=.35\textwidth]{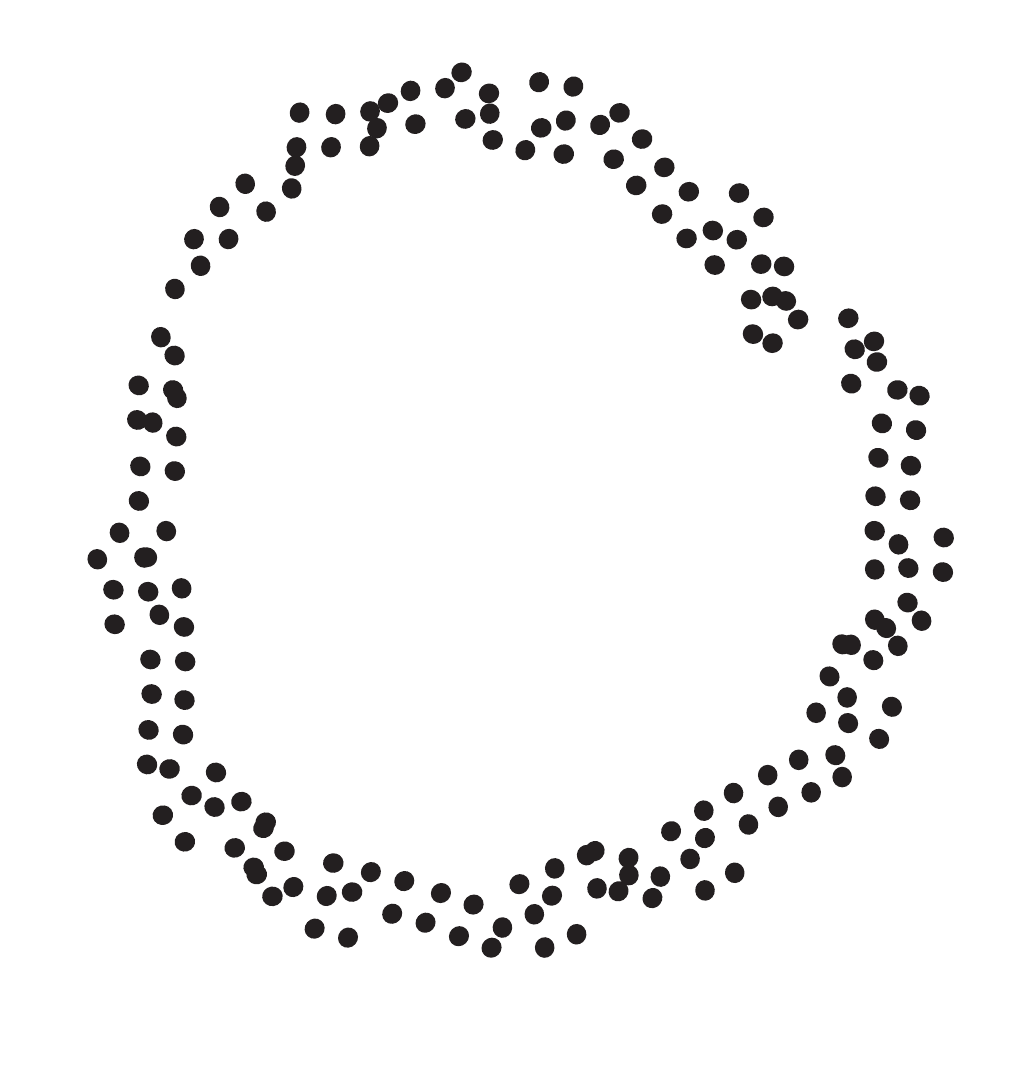}
\caption{Statistical Circle  }\label{statisticcircle}
 
\end{figure}}

 When we examine $X$, we observe that it appears to be sampled from a loop, and would like to have algebraic tools that capture the ``loopy" structure of the set.  Note that we are only given a discrete set of points, so direct application of homological constructions will only produce the homology of a finite set of points.  However, we could attempt to construct a space from the set  $X$, which in a sense fills in the gaps between the points.   We will need to use some additional information about the points, and that will in this case be restriction of the Euclidean  metric to  $X$.   One distance  based construction is the {\em Vietoris-Rips complex}. 
 
 \begin{Definition}\label{defvr} For any finite metric space $(X,d)$, and every $R \geq 0$, let ${\cal V}(X,R)$ denote the simplicial complex with vertex set equal to $X$, and such that $\{ x_0, \ldots , x_k\}$ spans a $k$-simplex if and only if $d(x_i,x_j) \leq R$ for all $0 \leq i<j\leq k$.  
 \end{Definition}
 Notice that if $R$ is  smaller than the smallest interpoint distance, then ${\cal V}(X,R)$ will be a discrete complex on the set $X$.  On the other hand, if $R$ is greater than the diameter of $X$, then ${\cal V}(X,R)$ is a full simplex on $X$.  For intermediate values, one obtains other complexes.  In this case, there is a range of values of $R$ in which ${\cal V}(X,R)$ has the homotopy type of a circle.  One could ask if there is a principled way to choose a threshold $R$ based on only the distances between the points.  After a great deal of experimentation, one finds that this is a very difficult, if not unsolvable problem.  A question that one can ask is if there is a more structured object that one can study which incorporates all the thresholds in a single object, which can be analyzed in a number of different ways. Statisticians have confronted this problem in an analogous situation.  
 
{\bf Hierarchical Clustering:} The {\em clustering problem}  in statistics is to determine ways to infer the set of connected components of a metric space $X$ from finite samples $Y$.  One of the approaches statisticians  developed was to compute $\pi _0{\cal V}(Y,R)$ for a choice of threshold $R$, but they confronted the analogous problem to the one we discussed above, namely the selection of  $R$.  They came to the conclusion that the problem of  selecting a threshold in a principled way is  not a  well posed problem, but managed to construct a structured output, called a {\em dendrogram}, which allowed them to study the behavior of all thresholds at once.  It is defined as follows.  For each threshold $R$, we obtain a set of components $\pi _0{\cal V}(Y,R) $, which yields a partition $\Pi _R$ of $Y$.  If we have $R \leq R^{\prime}$, then $\Pi _R$ is a refinement of $\Pi _{R^{\prime}}$.  One definition of a dendrogram structure on a finite set  $Y$ is a parametrized family $\{ \Pi _R \}_{R \geq 0}$ of partitions of $Y$, with the  property that $\Pi _R $ refines $\Pi _{R^{\prime}}$ whenever $R \leq R^{\prime}$, and so that for any partition $\Pi$, the intervals $\{R | \Pi _R = \Pi \}$ are either empty or  closed  on  the left and open  on  the right.  We assume that there is an $R_{\infty}$ so that $\Pi_{R_{\infty}}$ is the partition with one block, namely $Y$.  The reason for this terminology is that this information is equivalent to a tree $\frak{D}$ with a reference map to the non-negative real line. The tree $\frak{D}$ is defined  as follows.   The points of $\frak{D}$ are pairs $(c,R)$, where $c$ is a block in the partition $\Pi _R$, and $0 \leq R \leq R_{\infty}$.  We clearly have a reference map to $[0,R_{\infty}]$ given by $(c,R) \rightarrow R$. To define a topology on this set, we construct an auxiliary space $Z$, defined as 
$$ Z = \coprod_{c} [0,R_{\infty}]_c
$$
where $c$ is a block in the partition $\Pi_0$, and $[0,R_{\infty}]_c $ denotes a copy of the interval $[0,R_{\infty}]$ labelled by $c$.  There is a natural map $\varphi$ from $Z$ to $\frak{D}$ given by $(c,t) \rightarrow (\rho _t(c),t)$, where $\rho _t $ denotes the projection $Y/\Pi_0 \rightarrow Y/\Pi _t$. It is clear that $\varphi$ is a surjective map, and therefore that $\frak{D}$ is the quotient of $Z$ by an equivalence relation.  The topology on $\frak{D}$ is the quotient topology associated to the topology on $Z$.  It is easy to check that this topology makes $\frak{D}$ into a rooted tree. 
  A tree with a reference map can be laid out in the plane, as indicated below, and one can recover directly the clustering at any given value of $R$.  
  {\begin{figure}[!hbp]
\centering
\includegraphics[width=.5\textwidth]{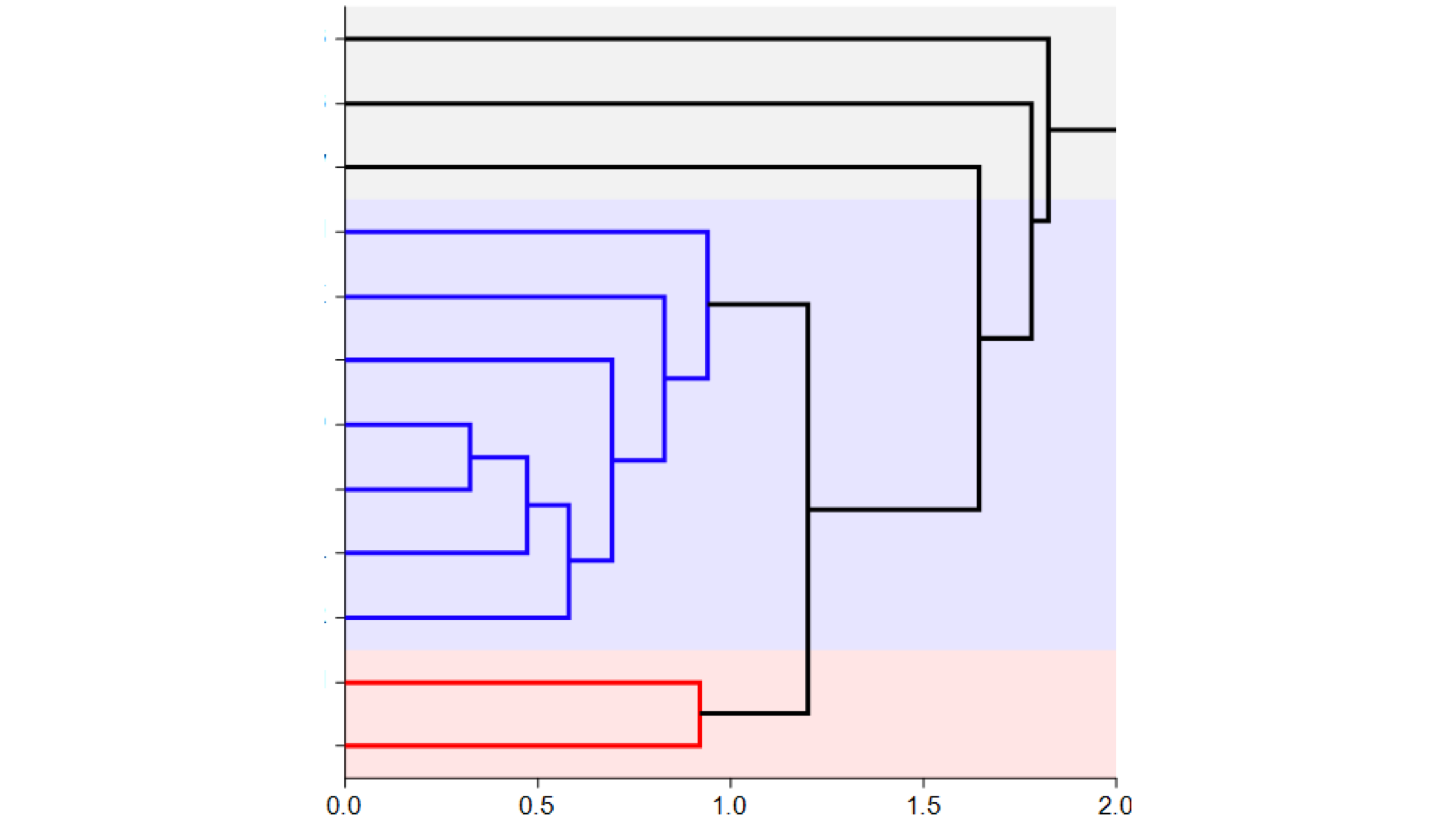}
\caption{Dendrogram}\label{dendrogramimage}
 
\end{figure}}

The dendrogram can be regarded as the ``right" version of the invariant $\pi _0$ in the statistical world of finite metric spaces.  The question now becomes if there are similar invariants that can capture the notions of higher homotopy groups or homology groups.  In order to define them, we need a preliminary definition. 

\begin{Definition} A submonoid $\Bbb{A} \subseteq \Bbb{R}_+$ is said to be {\em pure} if given any 
$r_0,r_1$, and $r_2$ in $\Bbb{R}_+$, with $r_0$ and $r_2$ in $\Bbb{A}$, so that $r_0 + r_1 = r_2$, then $r_1 \in \Bbb{A}$.  $\Bbb{A} $ is  a totally ordered set in its own right by restriction of the total order on $\Bbb{R}_+$.  
\end{Definition}
For example, $\Bbb{N}$ (the non-negative integers)  and $\Bbb{Q}_+$ are both pure.  We can now make a definition that includes the dendrogram as a special case.  

\begin{Definition} Let $\underline{C}$ denote any category. A persistence object in $\underline{C}$ is a functor $\Bbb{R}_+ \rightarrow \underline{C}$, where $\Bbb{R}_+$ denotes the ordered set of non-negative real numbers, regarded as a category, so that there is a unique morphism $r_0 \rightarrow r_1$ whenever $r_0 \leq r_1$. More generally, if $\Bbb{A} \subseteq \Bbb{R}_+$ is any {\em pure} submonoid, an $\Bbb{A}$-parametrized persistence object in $\underline{C}$ we will mean a functor from the  ordered set $\Bbb{A}$ to $\underline{C}$.  It is clear that the $\Bbb{A}$-parametrized  persistence objects in $\underline{C}$ form a category in their own right, where the morphisms are the natural transformations of functors.  We will denote this category by $\frak{P}\Bbb{A} \underline{C}$, with the special case of $\Bbb{R}_+$ denoted by $\frak{P}\underline{C}$.    \end{Definition} 

It is now easily checked that the functor $\Bbb{R}_+ \rightarrow \underline{\mbox{Sets}}$ defined by $$R \rightarrow \pi _0({\cal V}(X,R))$$ where $X$ is a finite metric space, is a persistence object in $\underline{\mbox{Sets}}$, or a {\em persistence set}.  This is the case because $\pi _0$ is a set valued functor.  Other topological invariants, such as homology or homotopy groups, take their values in the categories $\underline{\mbox{Grp}}$, $\underline{\mbox{Ab}}$, or $\underline{R-\mbox{mod}}$, and one can construct persistence objects in these categories by applying these functors to the Vietoris-Rips complexes.  The critical question becomes whether or not the isomorphism classes of these persistence objects are in any sense understandable, and useful for distinguishing or understanding the underlying metric spaces.

\section{The Structure Theorem for Persistence Vector Spaces over a Field} For the entirety of this section, $K$ will denote a field, which will be fixed throughout.   We are interested in the isomorphism classification of $\Bbb{A}$-parametrized persistence $K$-vector spaces, for pure submonoids $\bb{A}$ of $\bb{R}_+$.  This is very complicated in general, but is manageable for objects of $\frak{P}\underline{\mbox{Vect}}(K)$ satisfying a finiteness condition, which is always satisfied for the Vietoris-Rips complexes associated with finite metric spaces. 
Let $\Bbb{M}$ be any commutative monoid.  By an $\Bbb{M}$-graded $K$-vector space, we will mean a $K$-vector space $V$ equipped with a decomposition 
$$ V \cong \bigoplus _{\mu \in \Bbb{M}} V_{\mu}
$$
Given two $\Bbb{M}$-graded $K$-vector spaces $V_*$ and $W_*$, by their tensor product we will mean the tensor product  $V \column{\otimes}{K} W$ equipped with the $\Bbb{M}$-grading given by 
$$(V \column{\otimes}{K}W)_{\mu} = \bigoplus_{\mu  _1 + \mu _2 = \mu} V_{\mu _1} \otimes W_{\mu _2}
$$
and we write $V_* \column{\otimes}{K} W_* $ for this construction.  An $\Bbb{M}$-graded $K$-algebra is then an $\Bbb{M}$-vector space $R_*$ together with a homomorphism $R_* \column{\otimes}{K}R_* \rightarrow R_*$, satisfying the associativity and distributivity conditions.  An important example is the monoid $K$-algebra $K[\Bbb{M}]_*$, for which the grading is given by 
$$ K[\Bbb{M}]_{\mu} = K\cdot \mu
$$
It will be convenient to write $t^{\mu}$ for elements $\mu \in K[\Bbb{M}]_*$.  
 We define the notion of a $\Bbb{M}$-graded $R_*$-module $M_*$ in a similar way.  

We specialize to the situation $\Bbb{M} = \Bbb{A}$, where $\Bbb{A}$  is a pure submonoid $\Bbb{A} \subseteq \Bbb{R}_+$. We will demonstrate the classification of   $\Bbb{A}$-parametrized persistence modules  by using an equivalence of categories to the category   of $\Bbb{A}$-graded $K[\Bbb{A}]_*$-modules.

\begin{Proposition}\label{pclassify} Let $K[\Bbb{A}]_*$ denote the monoid algebra of $\Bbb{A}$ over $K$, for a pure submonoid $\Bbb{A} \subseteq \Bbb{R}_+$, regarded as a graded $K$-algebra. Let $\underline{G}(\Bbb{A},K)$ denote the category of $\Bbb{A}$-graded $K[\Bbb{A}]_*$-modules.  Then there is an equivalence of categories 
$$ \frak{P}\Bbb{A}\underline{ \mbox{\em Vect}}(K) \cong \underline{G}(\Bbb{A},K)
$$.  
\end{Proposition} 
\begin{Proof}  Given a functor $\theta: \Bbb{A} \rightarrow \underline{\mbox{Vect}}(K)$, we will denote by $M(\theta)$ the graded $K$-vector space 
$$ M(\theta) = \bigoplus _{\alpha \in A}  \theta (\alpha )
$$
where the elements of grading $\alpha$ are precisely the elements of the  summand $\theta (\alpha)$.  We now extend the vector space structure to a graded $K[\Bbb{A}]_*$-module structure by defining the action $t^{\alpha} \cdot: \theta (\alpha^{\prime}) \rightarrow \theta (\alpha + \alpha ^{\prime})$ to be equal to the linear transformation $$\theta(\alpha^{\prime} \rightarrow \alpha + \alpha ^{\prime}):\theta (\alpha ^{\prime}) \rightarrow \theta (\alpha + \alpha ^{\prime})$$
where $\alpha^{\prime}\rightarrow \alpha + \alpha ^{\prime}$ denotes the unique morphism in $\Bbb{A}$ from $\alpha $ to $\alpha + \alpha ^{\prime}$.  It is clear that this defines an $\Bbb{A}$-graded $K[\Bbb{A}]_*$-module structure on $M(\theta)$, and that $M$ is a functor.  We produce an inverse functor $\eta: \underline{G}(\Bbb{A}, K) \longrightarrow 
\frak{P}{\Bbb{A}}\underline{\mbox{Vect}}(K) $ on objects by setting 
$$
\eta(M_*) (\alpha) = M_{\alpha}
$$
and on morphisms  by 
$$
\eta(M_*)(\alpha^{\prime} \rightarrow \alpha + \alpha ^{\prime}) =t^{\alpha} \cdot:M_{\alpha ^{\prime}}\rightarrow M_{\alpha + \alpha ^{\prime}}
$$
The functors $\theta$ and $\eta$ are clearly inverse to each other. 
\end{Proof}

Let $\Bbb{A}$ be any pure submonoid of $\Bbb{R}_+$.  Then for any $\alpha \in \Bbb{A}$, we define $F(\alpha)$ to be the free $\Bbb{A}$-graded $K[\Bbb{A}]_*$-module on a single generator in grading $\alpha$.  For any pair $\alpha, \alpha ^{\prime} \in \Bbb{A}$, where $\alpha^{\prime} > \alpha$, we define $F(\alpha, \alpha ^{\prime})$ to be the quotient 
$$  F(\alpha)/ (t^{\alpha ^{\prime}-\alpha }\cdot F(\alpha ))
$$
The following result describes the isomorphism classification of finitely presented graded $K[\Bbb{A}]_*$-modules.  
\begin{Proposition}\label{classification} Any finitely presented object of $\underline{G}(\bb{A},K)$ is isomorphic to a module of the form 
$$ \bigoplus _{s=1}^{m} F(\alpha _s) \oplus \bigoplus _{t=1}^{n} F(\alpha _t,\alpha ^{\prime}_t)
$$
Moreover, the decomposition is unique up to reordering of summands. The kernel of any homomorphism between two finitely generated free $\Bbb{A}$-graded modules is itself a finitely generated free $\Bbb{A}$-graded module.  
\end{Proposition}
\begin{Remark} {\em This result is formally very similar to the structure theorem for finitely generated modules over a principal ideal domain (PID).  Indeed, for the case of $\Bbb{A} = \Bbb{N}$, where $K[\bb{A}]$ is Noetherian, the result is exactly a structure theorem for finitely generated graded modules over the graded ring $K[t]$.  For other choices of $\Bbb{A}$, $K[\Bbb{A}]_*$ is not necessarily  Noetherian.  However, it does turn out to be {\em coherent}, i.e. having the property that the kernel of any homomorphism between finitely generated free modules is finitely generated.  }
\end{Remark}

\begin{Proof} A proof is given in \cite{actanumerica}.  \end{Proof} 

{\bf Remark:} {\em Notice that the proof also gives an algorithm for producing a matrix in the diagonal form given above. }

The above theorem is summarized using the following definition.  

\begin{Definition} By an {\em $\bb{A}$-valued barcode}, we will mean a finite set of elements 
$$ (\alpha, \alpha \pr) \in \bb{A} \times (\bb{A} \cup \{ + \infty \})
$$
satisfying the condition $\alpha < \alpha \pr$.  An $\bb{A}$-valued barcode is said to be {\em finite} if the right hand endpoint $+ \infty$ does not occur.  If $\bb{A} = \bb{R}_+$, we will simply refer to it as a barcode, without labeling by the monoid.  We have shown that isomorphism classes of  finitely presented $K[\bb{A}]_*$-modules are in bijective correspondence with $\bb{A}$-valued barcodes.  
\end{Definition}  

\begin{Remark} {\em Barcodes are often represented visually as collections of intervals. The image in Figure \ref{barcodes}  shows barcodes in dimensions zero and one. 
{\begin{figure}[!hbp]
\centering
\includegraphics[width=.5\textwidth]{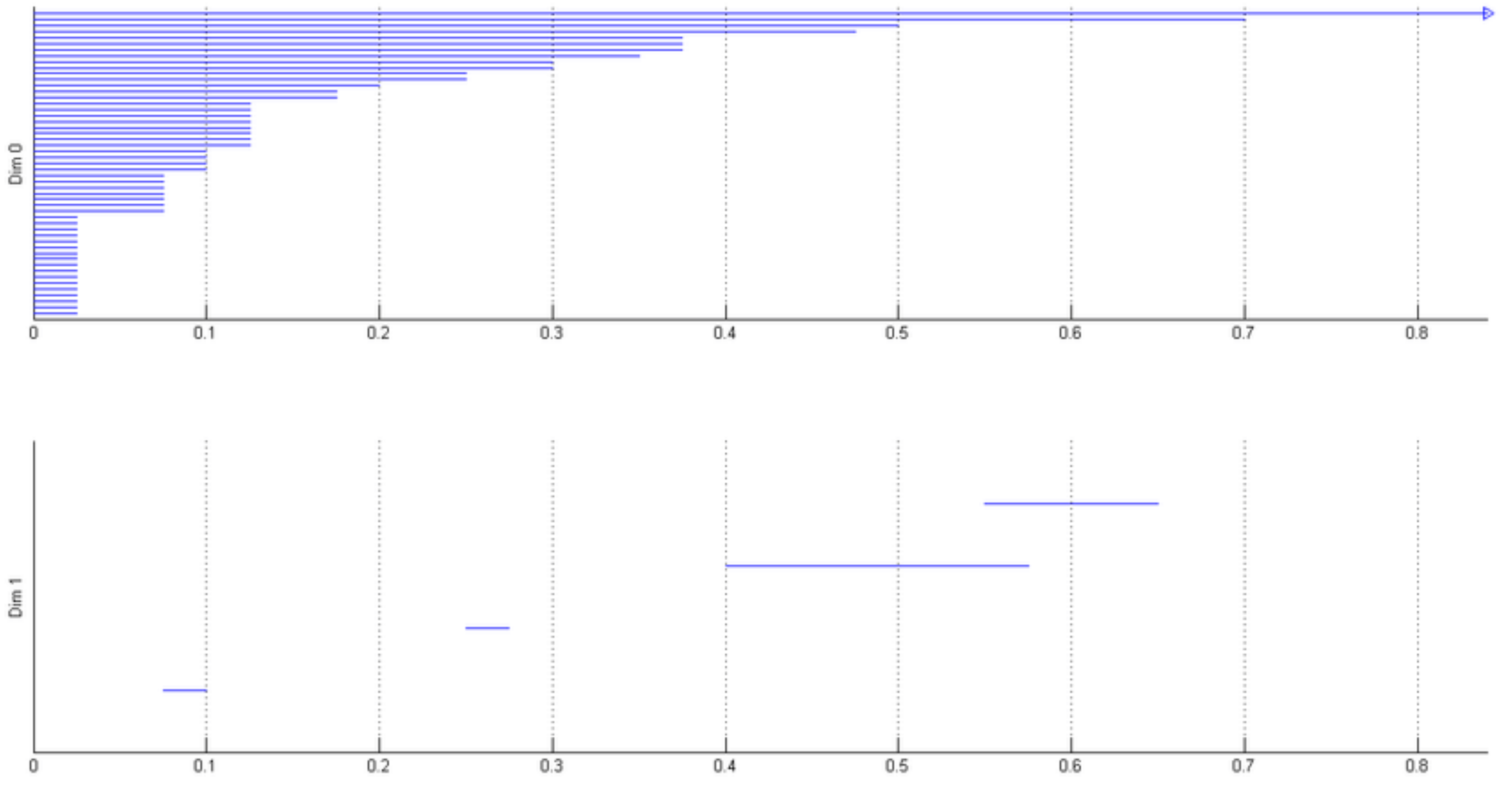}
\caption{Persistence Barcodes  }\label{barcodes}
 \end{figure}}

Notice that the zero-dimensional barcode has one infinite interval, while the one-dimensional barcode is finite.  There is an equivalent visual representation called the {\em  persistence diagram} in which each interval is encoded as a point $(x,y)$  in the plane, where $x$ and $y$ are the left and right hand endpoints respectively. An example  is  pictured  in Figure \ref{persistencediagram}. 
{\begin{figure}[!htp]
\centering
\includegraphics[width=.5\textwidth]{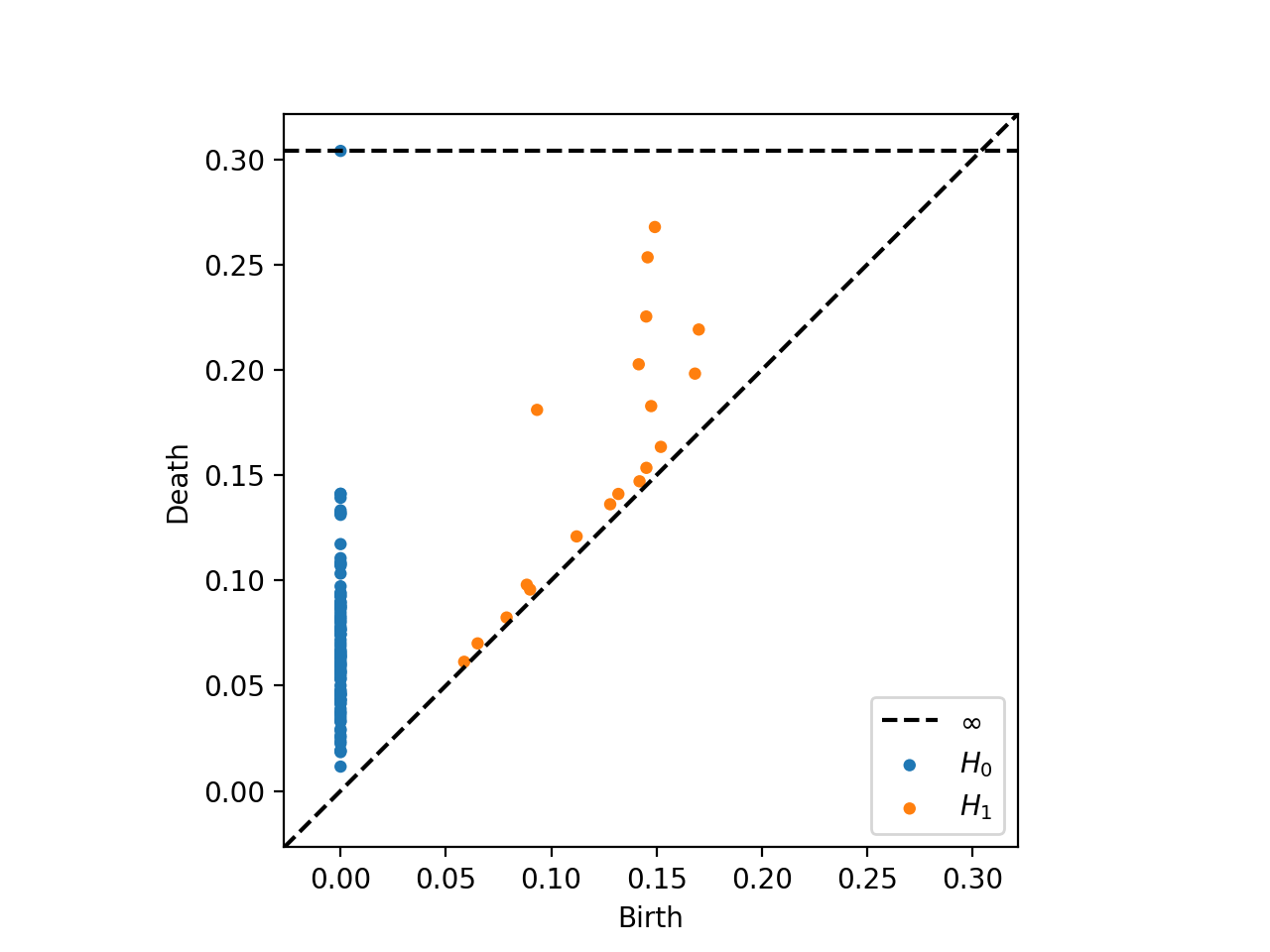}
\caption{Persistence Diagram}\label{persistencediagram}
 
\end{figure}}

An advantage of this representation made apparent in this image is that one can represent different homology groups in the same diagram.  The blue dots are in this case the zero-dimensional persistence diagram and the orange ones are the one-dimensional ones.  When there are infinite intervals, one often selects an upper threshold $\tau$ for the persistence parameter, and represents the infinite interval by one with left hand endpoint $\tau$.  }
\end{Remark}

We also have the following. 

\begin{Proposition} The kernel of a homomorphism between finitely generated free $K[\bb{A}]_*$-modules is finitely generated free.  
\end{Proposition}
\begin{Proof} This follows immediately from the matrix analysis in the proof of Proposition \ref{classification}.  
\end{Proof}
\begin{Corollary} \label{fp}  Given any chain complex of finitely generated free graded $K[\bb{A}]_*$-modules, the homology modules are finitely presented $K[\bb{A}]_*$-modules.  
\end{Corollary} 

We also interpret this result in terms of persistence vector spaces.  The category $\frak{P}{\Bbb{A}}\underline{\mbox{Vect}}(K)$ is clearly an abelian category.  For any $\alpha _0 \in \bb{A}$, we define $V(\alpha)$ to be the persistence vector space $\{ V _{\alpha} \}_{\alpha \in \bb{A}}$, where $V_{\alpha} = \{ 0 \}$ for $\alpha < \alpha _0$, $V_{\alpha } = K $ for $\alpha \geq \alpha _0$, and where for any $ \alpha _0 \leq \alpha \leq \alpha \pr$ the linear transformation $V_{\alpha} \rightarrow V_{\alpha \pr}$ is the identify on $K$.  For any $ \alpha < \alpha \pr$, we also define the persistence vector space $V(\alpha, \alpha \pr)$ to be the quotient of $V(\alpha)$ by the image of the natural  inclusion $V(\alpha \pr) \hookrightarrow V(\alpha)$.  Corollary \ref{fp} now has the following consequence. 

\begin{Corollary} Let $C_*$ be a chain complex of $\bb{A}$-persistence vector spaces, so that for every $s$, $C_s$ is a finite direct sum of persistence vector spaces, each of which is of the form $V(\alpha)$ for some $\alpha \in \bb{A}$.  Then for each $s$, $H_s(C_*)$ is isomorphic to a direct sum of finitely many persistence vector spaces, each of which is of the form $V(\alpha) $ or $V(\alpha,\alpha \pr)$ for $\alpha, \alpha \pr \in \bb{A}$, and $\alpha < \alpha \pr$.  In particular, the homology in each dimension $s$ is uniquely described by  an $\bb{A}$-valued barcode. 
\end{Corollary}

\section{Complex Constructors}
\subsection{Introduction} All data that we consider consists of finite sets of points $X$.   The space $X$ itself is uninteresting topologically, since it is a discrete set of points.   This means that we have to build a space using  auxiliary information attached to the set of points.  The auxiliary information we choose is a metric on the set $X$, so $X$ is a finite metric space.  In the context of data sets, metrics are often referred to as {\em dissimilarity measures}, since small distances between data points often reflect notions of similarity between data points.    Often the metric chosen is the restriction of a well known and analyzed  metric on an ambient space containing $X$, such as $n$-dimensional Euclidean space. Other choices that are often appropriate are Hamming distance, correlation distances, and normalized variants (mean centering of coordinate functions, normalizing variance to 1)  of Euclidean distance.   One method for constructing spaces based on metrics  is the Vietoris-Rips complex that we have seen  above.  It is actually a persistence object in the category of simplicial complexes.  All the constructions we will look at except the Mapper construction  are naturally persistence objects in the category of simplicial complexes.  The Mapper construction can also be equipped with many such structures, but they are not canonical.  Because of the presence of a persistence structure, the \v{C}ech, Vietoris-Rips, alpha, and witness complexes all have persistence barcodes associated to them in all non-negative dimensions.  It is obvious from the constructions that the zero-dimensional barcodes have a single infinite interval, and that all higher dimensional barcodes are finite.  
\subsection{\v{C}ech Construction}  Let $(X,d)$ denote a finite metric space.   Given a threshold parameter value $R$, let ${\cal U}_R$ denote the covering of $X$ by all balls $B_R(x) = \{ x^{\prime} | d(x,x^{\prime}) \leq R \}$.  The {\em \v{C}ech complex at scale $R$} is the nerve of the covering ${\cal U}_R$, and we denote it by $\check{C}(X;R)$.  It is clear that for $R \leq R^{\prime}$, there is an inclusion 
${\check{C}}(X;R) \hookrightarrow \check{C}(X,R^{\prime})$, and that therefore $\{\check{C}(X;R) \}_R$ is a persistence object in the category of simplicial complexes.  From the theoretical point of view, it has the advantage that given a covering of a topological space (with suitable point set hypotheses) by open sets, the nerve lemma (see \cite{kozlov}, Thm. 15.21) gives a criterion that guarantees that the nerve of the covering is homotopy equivalent to the original space.  
\subsection{Vietoris-Rips Complex}
The \v{C}ech construction from the previous section is computationally expensive, because it involves computing simplices individually at every level.  We can create another construction that has a strong relationship with the \v{C}ech construction.      Recall that a simplicial complex $Z$ is said to be a {\em flag complex} if for any collection $\sigma$ of vertices $\{ z_0, \ldots , z_k \}$ for which each pair $(z_i,z_j)$ is an edge, $\sigma$ is a $k$-simplex of $X$.  From the computational point of view, flag complexes are attractive  because one needs only enumerate all the edges in the complex, rather than all the  higher order simplices.  The Vietoris-Rips complex which was defined in Definition \ref{defvr} is by definition a flag complex for every parameter value $R$.  There is a  relationship between the persistent \v{C}ech and Vietoris-Rips complexes.  

\begin{Proposition} \label{interleave} There are  inclusions 
$$ \check{C}(X,R) \hookrightarrow {\cal V}(X,2R)
$$ 
and
$$ {\cal V}(X,R)  \hookrightarrow \check{C}(X,R)
$$
\end{Proposition}
For specific situations, this bound can be improved.  For example, it is shown in \cite{desilva1} that  if $X \subseteq \bb{R}^d$ is equipped with the restricted metric, then
$$ {\cal V}(X,{R ^{\prime}}) \subseteq \check{C}(X,R/2 ) \subseteq {\cal V}(X, \epsilon ) \mbox{ if } \frac{R}{R ^{\pr}} \geq \sqrt{\frac{2d}{d+1}}
$$
This result allows us to compare homology computed using the \v{C}ech and Vietoris-Rips methods.

\subsection{Alpha Complex}
There is another kind of complex whose dimension is low and which generally has a moderate number of simplices.  It is called the {\em alpha complex}, or the {\em alpha shapes complex}, and was introduced in \cite{alpha1}, with a thorough description in \cite{alpha}.  It applies to data $X$ that is embedded Euclidean space $\Bbb{R}^n$, and so that the metric on $X$ is the restriction of the Euclidean metric to $X$.  Typically the number $n$ is relatively small, say $\leq 5$, because the construction becomes quite  expensive in higher dimensions.  Also, the complex is generically of dimension $\leq n$.  The notion of {\em generic} is the following.  Given any set of points  $X\subseteq \Bbb{R}^n$, it is possible for the alpha complex to have dimension higher than $n$, but it is possible to perturb all the points by an arbitrarily small amount and obtain a complex that is $n$-dimensional.  

For any point $x \in X$, we define the {\em Voronoi cell} of $x$, denoted by $V(x)$, by 
$$  V(x) = \{ y \in Y| d(x,y) \leq d(x^{\prime} , y) \mbox{ for all } x^{\prime} \in X \}
$$
The collection of all Voronoi cells for a finite subset of Euclidean space is called its {\em Voronoi diagram}.  A Voronoi diagram in $\Bbb{R}^2$ might look like this.

{\begin{figure}[!htp]
\centering
\includegraphics[width=.3\textwidth]{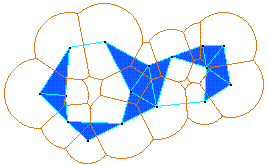}
\caption{Alpha Complex in $\Bbb{R}^2$}\label{voronoidiagram}
 
\end{figure}}

For each $x \in X$, we also denote by $B_{\epsilon}(x)$ the set $\{ y \in Y | d(x,y) \leq \epsilon \}$.  By the {\em $\alpha$-cell} of $x \in V(x)$ with scale parameter $\epsilon$, we will mean the set 
$A_{\epsilon}(x)  = B_{\epsilon}(x) \cap V(x)$.  The {\em $\alpha$-complex with scale parameter $\epsilon$}  of a subset $x \in X$, denoted by $\alpha _{\epsilon}(X)$ will be the abstract simplicial complex with vertex set $X$, and where the set $\{ x_0, \ldots , x_k \}$ spans a $k$-simplex if and only if 
$$  \bigcap _{i=0}^kA_{\epsilon} (x_i) \neq \emptyset
$$
It is of course the nerve of the  covering of $\bb{R}^n$ by the sets $A_{\epsilon}(x_i)$.  

\subsection{Witness Complex} 

The {\em witness complex} was introduced in \cite{witness}.  It  can be thought of as an analogue of the alpha complex for non-Euclidean data. The construction of the Voronoi cells is not dependent on the fact that the embedding space is Euclidean. It can be constructed for any embedding of a data set in a larger metric space.  Given a data set $X$, we therefore select a subset of {\em landmarks}  ${\cal L} \subseteq X$.  We can now build the analogues of the Voronoi cells for each of the landmark points within $X$, and construct the nerve of the covering.   In this case, both the ambient space and the landmark set are usually taken to be finite.  We will also need to introduce persistence into this picture.  The construction is as follows.  

\begin{Definition} Given a  metric space $(X,d)$,  a finite  subset ${\cal L} \subseteq X$, called the {\em landmark set}, and a persistence parameter $R$, and for every $x \in X$ we denote by $m_x$ the distance from $x$ to the set ${\cal L}$.  We define a simplicial complex $W(X,{\cal L},R)$ as follows.  The vertex set of $W(X,{\cal L},R)$ is the set ${\cal L}$, and  $\{ l_0, l_1, \ldots , l_k \}$ spans a $k$-simplex if and only if  there is a point $x \in X$ (the witness) so that $d(x, l_i)\leq m_x +R$ for all $i$.  The family of complexes $\{W(X,{\cal L},R)\}_R$ form a persistence simplicial complex.  \end{Definition}

   There are several  variants on this construction. For example, there is the ``lazy" version in which the $1$-simplices are the identical to the 1-simplices of the witness construction, but in which we declare that any higher dimensional simplex is an element of the complex if and only if all its one dimensional faces are.  Each of the lazy complexes is a flag complex. The lazy witness complex bears the same relationship to the witness complex as the Vietoris-Rips complex bears to the \v{C}ech complex.     There is also the {\em weak witness complex}, $ \{W^{weak}(X,{\cal L},R) \}_R$, which is defined as follows.  For each point $x \in X$, we let $\delta _x$ denote the  distance to the {\em  second } closest element of $\Lambda$ to $x$.  We then declare that a pair $(\lambda _1, \lambda _2)$ spans an edge of $W^{weak} (X,{\cal L},R)$ if and only if there is an $x \in X$ so that

$$\max(d(\lambda _1,x), d(\lambda _2,x)) \leq \delta _x + R
$$

A higher dimensional simplex $\{\lambda _0, \ldots , \lambda _n \}$ is contained in $W^{weak} (X,{\cal L},R)$ if an only if all of its edges are        
contained in it.  This is a very useful construction because the persistence ``starts faster" than the standard complex.  It is often the case that one obtains the ultimate result even at $R =0$, or for very small values of $R$.
\subsection{Mapper} 
Another construction is based on Morse theoretic ideas. We motivate it by considering a space level construction.  Suppose that we have a continuous map $r:X \rightarrow B$ of spaces, and suppose further that $B$ is equipped with an open covering ${\cal U} = \{ U_{\alpha } \}_{\alpha \in A}$.  We obtain the open covering $r^{-1}{\cal U} = \{ r^{-1}U_{\alpha} \}_{\alpha \in A}$ which can be refined into a new  covering $r^{-1}{\cal U}^*$  by decomposing each set $r^{-1}U_{\alpha }$ into its connected components.  We note that the dimension of the nerve of $r^{-1}{\cal U}^*$ is less than or equal to the dimension of the nerve of ${\cal U}$, so this construction, like the alpha complex, produces complexes of bounded dimension.   The analogue of this construction for finite metric spaces is obtained by assuming that the finite metric space is equipped with a map $\rho$  to a reference space $B$, and replacing the connected component construction by the output of a clustering algorithm.  A simple algorithm  to use is single linkage  hierarchical clustering, where one makes a choice of threshold based on a simple heuristic, such as the one found in \cite{mapper}.  Once this is done, one obtains a covering of the finite set $X$ by the collection of all clusters constructed in each of the sets $\rho ^{-1}(U_{\alpha})$, and constructs the nerve complex.  This construction is referred to as {\em Mapper}. Usually, $B$ is chosen to be $\bb{R}^n$, for a small positive integer $n$, and therefore the reference map is determined by an $n$-tuple of real valued functions on the metric space $X$.  The reference maps can be chosen  in many ways, giving different views of the data.  Some standard choices are density estimators, measures of centrality, coordinates in linear algebraic algorithms such as principal component analysis and multidimensional scaling (\cite{friedman}), or individual coordinates when a metric space is obtained as a subspace of $\bb{R}^N$ for some $N$.  The method has been used extensively in work on life sciences data sets, see for example \cite{dudley}, \cite{nicolau}, \cite{olin},  and \cite{saggar}.

\section{Metrics on Barcode Space, and  Stability Theorems}\label{wholemetrics}
Since persistent homology is used to analyze data sets, and data sets are often noisy in the sense that one does not want to assign meaning to small changes, it is important to analyze the stability of persistent homology outputs to small changes in the underlying data.  In order to do this, it is very useful to construct metrics on the output barcodes so that one can assert the continuity of the assignment of a barcode to a finite metric space or to the graph of a function.  Informally, one wants to prove that small changes in the data give rise to small changes in the barcodes.  Since small changes will often result in a change in the number of bars, it will be important to construct a set of all barcodes, on which we can impose a metric.  The following is a natural construction.    Let $n$ be a positive integer, and let $B_n$ denote the set of unordered  $n$-tuples of closed intervals $[x,y]$, where we permit $x = y$, and require $x,y \geq 0$.  It is understood that $B_0$ consists of a single point, namely the empty set of intervals.   The set $B_n$  can be described as the orbit space of the $\Sigma _n$-action on the set ${\cal I}^n$ which permutes coordinates, where ${\cal I}$ denotes the set of closed intervals.  To consider all barcodes, we form 
$$ \frak{B}_+ =  \coprod _{n \geq 0}B_n
$$
and define the full barcode space $\frak{B}$ to be the quotient $\frak{B}_+/\sim$, where $\sim$ is the equivalence relation generated by the relations 
$$ \{I_1, \ldots,I_{k-1}, [x_i,x_i] ,I_{k+1}, \ldots ,I_n\} \sim 
 \{ I_1, \ldots,I_{k-1},I_{k+1} , \ldots ,I_n \} 
$$
The idea is that intervals of length zero, which do not represent non-zero vector spaces, should be ignored.  We will need to construct metrics on $\frak{B}$, and to prove continuity results for these metrics.    
\subsection{Metrics on Barcode Space} \label{metrics} 
The general idea for the construction of metrics on barcode space is to consider the set of all partial matchings between the intervals in the barcodes, assign a penalty to each such matching, and finally to minimize this penalty over the set of all matchings.  Partial matchings are a little awkward, so for a pair of barcodes $B_1$ and $B_2$  one instead considers actual bijections $ B_1^{\pr} = B_1 \cup Z \rightarrow B_2 \cup Z = B_2^{\pr}$, where $Z$ is the set consisting of a countable number of copies of the interval of length zero $[x,x]$ for each $x \geq 0$.  To start, one  assigns a penalty  $\pi([x_1,x_2],[y_1,y_2]) = ||(y_1 - x_1, y_2 -x_2)||_{\infty}$ for every pair of intervals including those of length zero.   Next, given two barcodes $B_1$ and $B_2$, let ${\cal{D}}(B_1, B_2)$ denote the set of all bijections $\theta: B_1^{\pr} \rightarrow B_2^{\pr}$ for which $\pi(I,\theta (I)) \neq 0$ for only finitely many $I \in B_1^{\pr}$.  Given a positive number $p$, we  extend the penalty function $\pi$  from individual intervals to barcodes by forming 
$$ W_p(B_1,B_2) = \mathop{\mbox{inf}}_{\theta \in {{\cal{D}}(B_1,B_2)}} (\sum _{I \in B_1^{\pr}}\pi (I, \theta (I))^p)^{\frac{1}{p}}
$$
As usual, $p= \infty$ is interpreted as the $L_{\infty}$ norm. 
\begin{Definition} For any $p >0$ and also for  $p = \infty$, we refer to $W_p(B_1,B_2)$ as the {\em $p$-Wasserstein distance}.  For $p = \infty$, this distance is often referred to as the {\em bottleneck distance}.  
\end{Definition}  
It is readily verified that under this definition, $W_p$ defines a metric on $\frak{B}$.  
\subsection{Stability Theorems}\label{stable}
The metrics defined in the previous section have stability theorems associated to them, that assert that the assignment of barcodes is a continuous process.  There are two  situations of interest.  
\begin{enumerate}
\item{Gromov has defined (see \cite{gromov}) a metric $d_{GH}$ on the set of all isometry classes of compact metric spaces, called the Gromov-Hausdorff metric.  For any integer $k \geq 0$, one can view the assignment to any finite metric space its barcode as a map from the set of isometry classes of finite metric spaces to $ \frak{B}$, and one can ask about its continuity properties.  }
\item{Let $f:X \rightarrow \bb{R}_+$ denote a continuous function.  One can assign to each such $f$ and each integer $i \geq 0$ the persistence vector space $\{ H_i(f^{-1}([0,r]))\}_{r \geq 0}$.  Under suitable situations (e.g. where $X$ is a finite simplicial complex and $f$ is linear on simplices, or where $M$ is a compact manifold and $f$ is smooth), one can show that the associated barcode is of finite type.  It is then an interesting question to ask about the continuity properties of this assignment, where one assigns various metrics to the set of functions on $X$.  }
\end{enumerate} 
There are theorems in both these cases.  The following theorem demonstrates the continuity of the assignment of a $k$-dimensional barcode with coefficients in a field to a finite metric space, when the metric on the set of isometry classes of finite metric spaces is the Gromov-Hausdorff distance.  
\begin{Theorem}{\em \cite{ghstable}} For any two finite metric spaces $X$ and $Y$ and integer $k \geq 0$, let $B(X)$ and $B(Y)$ denote the $k$-dimensional barcodes for the Vietoris-Rips complexes of  $X$ and $Y$ in a field $K$.  Then 
$$ W_{\infty}(B(X),B(Y)) \leq d_{GH}(X,Y)$$
\end{Theorem} 

There is a direct analogue for continuous real valued  functions on a topological space.  In order to state it, we need a pair of definitions.  

\begin{Definition} Let $X$ be a topological space and $f$ a continuous real valued function on $X$.  A real number $a$ is said to be a {\em homological critical value} of $f$  if for some $k$ and  all sufficiently small $\varepsilon >0$ the inclusion 
$$ H_k(f^{-1}(- \infty ,a-\varepsilon]) \rightarrow H_k(f^{-1}( - \infty , a + \varepsilon ])
$$ 
is not an isomorphism.  
\end{Definition} 
\begin{Definition}  Let $X$ be a topological space and $f$ be a continuous real-valued function.  We say $f$ is {\em tame} if there are finitely many homological critical values and all the homology groups $H_k(f^{-1}(- \infty, a]))$ with coefficients in a field $K$ are finite dimensional.  
\end{Definition} 
\begin{Remark}{\em Tameness holds in  a number of familiar situations. 
\begin{enumerate}
\item{Morse functions on compact smooth manifolds.}
\item{Functions on finite simplicial complexes that are linear on each simplex.  }
\item{Morse functions on compact Whitney-stratified spaces.}
\end{enumerate}   }
\end{Remark}
The theorem is as follows.  
\begin{Theorem}\label{bstable}{\em \cite{bottlestable}} Let $X$ be triangulable space, and suppose $f,g: X \rightarrow \bb{R}$ are tame continuous functions.  Let $B(f)$ and $B(g)$ denote the barcodes attached to  $\{ H_k(f^{-1}(- \infty , r];K ))\}_{r}$ and $\{ H_k(g^{-1}(- \infty , r];K ))\}_{r}$, respectively.  Then 
$$ W_{\infty}(B(f),B(g)) \leq ||f-g||_{\infty}
$$
\end{Theorem} 

The situation for the $p$-Wasserstein distances where $p \neq \infty$ is more complex.  We will need constraints on the metric space as well as on the functions.  For the space $X$, it is required to be a triangulable compact metric space, so $X$ is homeomorphic to a finite simplicial complex.  In addition, though, there is a requirement that the number of simplices required to construct a triangulation where the diameter of the simplices is less than a threshold $r$.   Specifically, for a given $r >0$, we define $N(r)$ to be the minimal number of simplices in a triangulation of $X$ for which each simplex has diameter $\leq r$.  We will assume that $N(r)$ grows polynomially with $r^{-1}$, i.e. that there are constants $C$ and $m$ so that $N(r) \leq \frac{C}{r^m}$.  It is easy to observe that this result holds for a finite simplicial complex $X$ equipped with the Euclidean metric obtained by restricting along a piecewise linear embedding $X \hookrightarrow \bb{R}^n$, as well as for a compact Riemannian manifold.  In \cite{lipschitz}, it is proved that any metric space satisfying this condition also satisfies a homological condition.  To state the homological condition, given a barcode $B = \{[x_1,y_1], \ldots , [x_n,y_n] \}$, we define $P_k(B)$ to be the sum 
$$ P_k(B) = \sum _i (y_i - x_i )^k
$$

\begin{Lemma}\label{condition} {\em \cite{lipschitz}} Suppose that $X$ is as above, and that $k$ is a nonnegative real number.  Then there is a constant $C_X$ so that $P_k(B(f)) \leq C_X$ for every tame function $f : X \rightarrow \bb{R}$ with Lipschitz constant $L(f) \leq 1$, where $B(f)$ is defined as in Theorem \ref{bstable}.  
\end{Lemma} 
\begin{Definition} When the conclusion of Lemma \ref{condition} holds for a metric space $X$ and real number  $k\geq1$, we say that $X$ {\em implies bounded degree-$k$ persistence}. 
\end{Definition}The theorem is now as follows.  
\begin{Theorem}{\em \cite{lipschitz}}  Let $X$ be a triangulable metric space that implies bounded degree-$k$ persistence for some real number $k \geq 1$, and let $f,g:X \rightarrow \bb{R}$ be two tame Lipschitz functions.  Let $C_X$ be the constant in Definition \ref{condition}, and let $L(f)$ and $L(g)$ denote the Lipschitz constants for $f$ and $g$ respectively.   Then
$$ W_p(f,g) \leq C^{\frac{1}{p}}\cdot ||f-g||_{\infty}^{1- \frac{k}{p}}
$$
for all $p \geq k$, where $C = C_X\mbox{max} \{ L(f)^k, L(g)^k \}$. \end{Theorem}

\section{Tree-like Metric Spaces} Persistent homology gives ways of assessing the shape of a finite metric space.  One situation where this is very useful is in problems in evolution.  The notion that there is a ``tree of life" is a very old one which actually predates Darwin.  Different organisms of the same type have attached to them sequences of the same length  in a genetic alphabet $A$.  Therefore any set of organisms produces a subset of the set of sequences of fixed length in $A$.    One can assign a metric to the space of all such sequences using Hamming distances or variants thereof.  The notion that there is a tree of the various organisms within a fixed type can be restated in mathematical terms as stating that the space ${\cal S}(A)$ corresponding to the organisms in the family is well modeled by a {\em tree-like metric space}, i.e. a metric space which is isometric to the set of nodes in a tree, possibly with weighted edges,  equipped with the distance function that assigns to each pair of vertices of the tree the length of the shortest edge path between them.  This approximability could be called the {\em phylogenetic hypothesis} for the particular class of organisms.  Testing this hypothesis for particular genetic data sets  has usually been done by attempting to fit trees to a given metric space and attempting to assess how well the approximation fits.  Given the persistent homology construction, one is tempted to develop criteria attached to the barcodes that can distinguish between tree-like and non-tree-like metric spaces.  The following theorem, proved in \cite{evolution} gives such a criterion. 

\begin{Theorem} Let $X$ be  a finite tree-like metric space.  Then the $k$-dimensional persistent homology of $X$  vanishes for $k>0$.  
\end{Theorem}
\begin{Remark} {\em This theorem was proved in the context of a study of data sets of viral sequences.  In that paper it was also shown that representative cycles for generators of persistent homology  in positive degrees gave important clues to the mechanism of the failure of the  phylogenetic hypothesis.   }
\end{Remark} 
\section{Persistence and Feature Generation}\label{coordinatization}
\subsection{Introduction} \label{featureintro} 
The output of persistent homology is an interesting data type, consisting as it does of finite collections of intervals.  When humans are directly interpreting barcodes, they are typically able to interpret barcodes directly from this description.  However, there is a whole class of problems where computers are used to ``interpret" the barcodes.  For example, suppose that we have a database of complex molecules.  Each molecule is given as a collection of atoms and bonds, and the bonds may be equipped with lengths.  The set of atoms in a molecule  can be endowed with a metric by considering the edge-path distance using the lengths of the bonds as the lengths of the edges.  What we have is now a data set in which each of the {\em data points} is a finite metric space, and therefore possesses a barcode.  If there are many molecules, we cannot hope to interpret these barcodes by eye, and must therefore allow a computer to deal with them.  Machine learning algorithms are generally not well equipped to deal with data points described as sets, and it is therefore important to encode them somehow as vectors, which are the natural input to such algorithms.  In this section we will describe three distinct methods for ``vectorizing" barcodes, i.e. for creating coordinates on the set of barcodes. Specifically, we will define coordinates on the space $\frak{B}$ constructed in Section \ref{wholemetrics}.

\subsection{Algebraic Functions} \label{algebraic} This method proceeds from the observation that the sets $B_n$ can be viewed as subsets of a real algebraic variety.  The set ${\cal I}$ embeds as a subset of the two-dimensional affine space $\bb{A}^2(\bb{R})$, and is defined by the inequalities $x,y \geq 0$ and $y \geq x$ for an interval coordinatized by the pair $(x,y)$.  Consequently, we have an embedding 
$${\cal I}^n \hookrightarrow  \bb{A}^2(\bb{R})^n \cong \bb{A}^{2n}(\bb{R})
$$
and it is equivariant with respect to the permutation actions on ${\cal I}^n$ and $\bb{A}^2(\bb{R})^n$.  Under the identification $\bb{A}^2(\bb{R})^n \cong \bb{A}^{2n}(\bb{R})$, with $\bb{A}^{2n}(\bb{R})$ coordinatized using coordinates $(x_1, \ldots , x_n , y_1, \ldots , y_n )$, the corresponding action simply permutes the $x_i$'s and $y_i$'s among themselves.  It is a standard result in algebraic geometry (see \cite{invariant}) that for any action of a finite group on an affine algebraic variety (over $\bb{R}$ in this case), there is an orbit variety, whose affine coordinate ring is the ring of invariants of the group action. It is easy to verify that in this case, the orbit set of the action on the closed real points of $\bb{A}^{2n}(\bb{R})$ is exactly the symmetric product $Sp^n(\bb{R}^2)$.  Since $B_n \subseteq Sp^n(\bb{R}^2)$, elements in the affine coordinate ring of the orbit variety can be regarded as functions on $B_n$, so we now have an algebra of functions ${\cal A}_n$ on $B_n$.  This means that we can describe functions on the sets of barcodes with exactly $n$ intervals.  What one really wants is a ring of functions on all of $\frak{B}$. In order to construct such a ring, we observe that $\frak{B}$ can be described as a quotient of the  direct limit of the  system 
\begin{equation} \label{directsystem} B_0 \rightarrow B_1 \rightarrow B_2 \rightarrow \cdots 
\end{equation}
where the inclusion $B_n \rightarrow B_{n+1}$ is given by 
$$ \{[x_1,y_1],[x_2,y_2],\ldots ,[x_n,y_n] \} \rightarrow \{[x_1,y_1],[x_2,y_2],\ldots ,[x_n,y_n], [0,0] \} 
$$
There is a corresponding direct system of affine schemes 
\begin{equation} \label{schemedirect} Spec({\cal A}_0) \rightarrow Spec({\cal A}_1) \rightarrow Spec({\cal A}_2) \rightarrow \cdots 
\end{equation} 
which is compatible with the system (\ref{directsystem}) above.  The colimit of the system (\ref{schemedirect}) is an affine scheme, whose affine coordinate ring is the inverse limit of the system 
$$ {\cal A}_0 \longleftarrow {\cal A}_1 \longleftarrow {\cal A}_2 \longleftarrow \cdots 
$$
which we will denote $\overline{{\cal A}}$.  The ring $\overline{{\cal A}}$ can be analyzed, but is a bit too complicated to be used in applications.  To define a smaller subring, we note that $Spec(\overline{\cal A})$ is equipped with an action by the algebraic group $\bb{G}_m$, and we can define a subring $\overline{\cal A}^{fin}$  to consist of all those functions $f$ so that all the translates of $f$ under the $\bb{G}_m$-action span a finite dimensional vector subspace within $\overline{\cal A}$.  The ring $\overline{\cal A}^{fin} $ is actually a graded ring, since  the $\bb{G}_m$-action determines a grading on it.  Within $\overline{\cal A}^{fin}$ we define a subring ${\cal A}^{fin}$ which consists of all elements of 
$\overline{\cal A}^{fin}$ that respect the equivalence relation $\sim$.  The main result of \cite{adcock} is the following. 

\begin{Theorem} The rings in question have the following properties. 
\begin{enumerate}
\item{The ring $\overline{\cal A}^{fin}$ has the structure 
$$ \overline{\cal A}^{fin} \cong \bb{R}[x_{i,j}; 0 \leq i, 0 \leq j,  \mbox{ and } i+j>0]
$$}
\item{The subring ${\cal A}^{fin}$ is identified with 
$$ \bb{R}[x_{i,j}; 0 < i, 0 \leq j,  \mbox{ and } i+j>0]
$$}

\item{The element $x_{i,j}$ is the function given on a barcode $\{[x_1,y_1], \ldots , [x_n,y_n] \}$ by 
$$   \sum_{s=1}^n  (y_s-x_s)^i(y_s + x_s)^j
$$}
\item{The elements of ${\cal A}^{fin}$ separate points in $\frak{B}$}
\item{The ring ${\cal A}^{fin}$ injects into the ring of functions on $\frak{B}$.  }
\end{enumerate} 
\end{Theorem} 

We remark that these functions are not continuous for the bottleneck distance on $\frak{B}$.  In \cite{kalisnik}, a tropical version of this work is presented, which gives functions which are continuous for the bottleneck distance.  For the $p< \infty$ situation, the functions are continuous for the  $p$-Wasserstein distance if one gives assigns $\frak{B}$ the direct limit topology associated to the filtration of $\frak{B}$ by the images of the spaces $B_n$, defined in Section \ref{metrics}. 
  
   \subsection{Persistence Landscapes}
   {\em Persistence landscapes} were introduced in \cite{bubenik} as another vectorization of barcodes.  The vectorization in this case consists of an embedding of the set $\frak{B}$ in a set of  sequences of functions on the real line.  Let $(a,b)$ denote a pair of elements of 
   $ \bb{R} $
  with $a \leq b$.  Then we define a function $f_{a,b}(t)$ on the real line by setting 
   $$ f_{(a,b)}(t) = \mbox{min}(t-a,b-t)_+
   $$
   where $c_+ = \mbox{max}(c,0)$.  A quick analysis of $f_{(a,b)}$ shows that it is zero  for $t \leq a$ and $t \geq b$, that on the interval $[a,\frac{a+b}{2}]$ it is equal to the graph of a line of slope $1$ including the point $(a,0)$, and on the interval $[\frac{a+b}{2},b]$ it is the graph of a line of slope $-1$ including the point $(b,0)$.  The shape of the graph is that of a pyramid.
   
   {\begin{figure}[!htp]
\centering
\includegraphics[width=.5\textwidth]{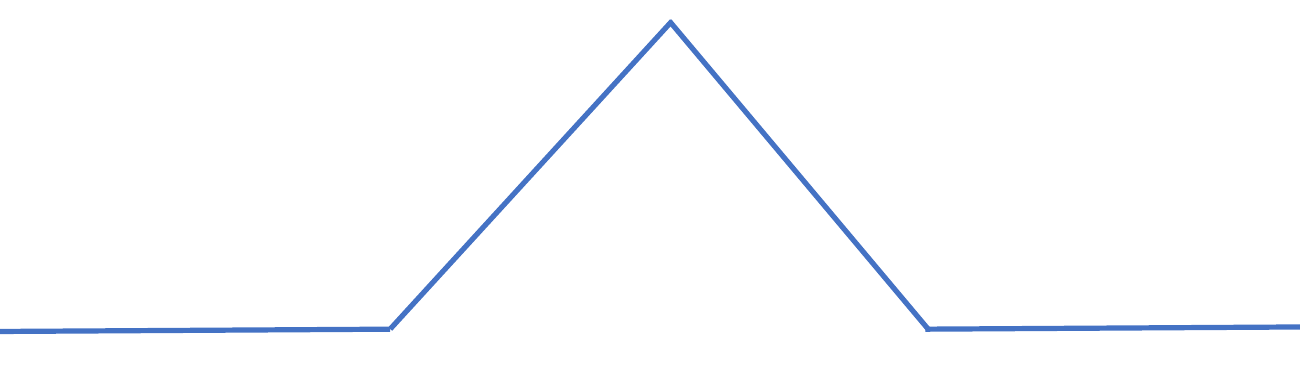}
\caption{Graph of $f_{a,b}$ }\label{pyramid}
 \end{figure}}

\begin{Remark}\label{compatible}  {\em Note that for $a = b$, $f_{(a,b)}(t) \equiv 0$. }
\end{Remark} 

Given a persistence barcode $B = \{(a_1,b_1), \ldots , (a_n,b_n) \}$, we define a family of functions $\lambda _k(t)$ parametrized by a positive integer $k$.  For $k = 1$, we let $\lambda _1 (B)  (t)$ denote the maximum of all the values $f_{(a_i,b_i)}(t)$ over all $i$.  For $k > 1$, we set $\lambda _k(B)  (t)$ equal to the $k$-th largest value occurring in the set $\{ f_{(a_i,b_i)}(t) \}_i$.  The family of functions $\{ \lambda _k (B) (t) \} _{k>0}$ is defined to be the persistence landscape of the barcode $B$.  It is clear that 
$$ \lambda _k(B)(t) \geq 0$$
that 
$$ \lambda _k(B)(t) \geq \lambda _{k + 1}(B)(t)
$$
and that 
$$ \lambda _k (B)(l) = 0 \mbox{ for } k > n
$$
In \cite{bubenik} it is also proved that each function $\lambda _k(B)$ is 1-Lipschitz, i.e. that 
$$ |\lambda _k(B)(t) - \lambda _kB) (t^{\prime})| \leq | t - t^{\prime}|
$$
To summarize, the persistence landscape lies in the vector space $\frak{F}$ of real valued functions on $\bb{N} \times \bb{R}$, and it follows directly from Remark \ref{compatible} that the definition gives us a function $PL:\frak{B} \rightarrow \frak{F}$. 

One extremely useful fact about the persistence landscape is that it is compatible with the various distances assigned to barcode spaces.  We let $\frak{F}_p\subset \frak{F}  $ denote the space of functions with finite $L_p$-norm $|| \mbox{ } ||_p$.   It is clear that the function $PL$ takes values in $\frak{F}_p$ for all $p$.    Recall the definition of the $p$-Wasserstein distance $W_p$  between barcodes from Section \ref{metrics}.  Bubenik now proves the following in \cite{bubenik}. 

\begin{Theorem} The function $d_p(B,B ^{\prime}) = ||PL(B) - PL(B^{\prime})||_p$ is a metric on $\frak{B}$.  The two metrics $W_{p+1}$ and $d_p$ generate the same topology on $\frak{B}$.  It follows that  the  map $PL$ is continuous when $\frak{B}$ is equipped with the metric $W_{p+1}$ and $\frak{F}_p$ is equipped with the metric associated with the $L_p$ norm. 
\end{Theorem} 

\begin{Remark}{\em Bubenik also provides explicit inequalities involving the two metrics in \cite{bubenik}}
\end{Remark} 
Bubenik also proves the following continuity  theorem.
\begin{Theorem} The map $PL$ is 1-Lipschitz from $\frak{B}$ equipped with the bottleneck distance to the space of persistent landscapes equipped with the sup norm distance.  This is equivalent to the algebraic statement
$$ |\lambda _k(B)(t) - \lambda _k(B \pr) (t)| \leq W_{\infty}(B, B \pr )
$$
\end{Theorem} 
\begin{Remark}{\em The map $PL$ separates points.    }
\end{Remark} 


\subsection{Persistence Images} There is another approach  that proceeds by treating a barcode, recoded as a persistence diagram, as a collection of point masses and then smoothing the corresponding measure to produce and image, which is finally discretized by selecting pixels and assigning each pixel the average value of the function on a box surrounding that pixel.  It is described in \cite{images}.  The detailed description is given in several steps. The input is a persistence diagram (it is more natural to use the persistence diagram view in this case), a collection of points $B = \{ (x_1,y_1), \ldots , (x_n,y_n)\}\subset \bb{R}^2$.  We assume we are given  a function $\phi: \bb{R}^2 \times \bb{R}^2 \rightarrow \bb{R}$ so that (a) $\phi _u = \phi (u, -)$ is a probability distribution on $\bb{R}^2$ for each $u \in \bb{R}^2$ and (b) the mean of $\phi _u$ is $u$. A standard choice  is that of a spherically symmetric Gaussian with mean $u$ and a fixed variance $\sigma$. We also assume we are given a continuous and piecewise differentiable nonnegative weighting function $f: \bb{R}^2 \rightarrow \bb{R}$ that  is zero along the $x$-axis.  

\begin{itemize}
\item{Apply the coordinate change $(x,y) \rightarrow (x,y-x) = (\xi, \eta)  $ to $\bb{R}^2$, to obtain the  new set of points 
${\cal B} = \{ (\xi _1, \eta _1)\ldots, (\xi _n, \eta _n)\}$ 
of the same cardinality in $\bb{R}^2$.  The points which correspond to short intervals are now all located near the $\xi $-axis.  The $\xi $-axis itself corresponds to intervals of length zero.  }
\item{Construct a new function $\rho_{\cal{B}}(z): \bb{R}^2 \rightarrow \bb{R}$, called the {\em persistence surface of ${\cal{B}}$}, defined by 
$$ \rho (z) = \sum _{i=1}^n f(\xi _i, \eta _i) \phi ((\xi_i, \eta _i), z) 
$$
Notice that $\rho $ vanishes on the $x$-axis.  }
\item{To construct a finite dimensional representation, we first assume that the persistence diagrams we will be dealing with will always lie in a bounded region in $\bb{R}^2$, and divide a box containing that region into a square grid.   Construct the vector with coordinates  in one to one correspondence with the squares of the grid, and assign  the entry corresponding to a square to be  the integral over that  square of $\rho_{\cal B}$. }
\end{itemize} 
The following is proved in \cite{images}
\begin{Theorem} The map $PI: \frak{B} \rightarrow \bb{R}^N$ which assigns to a persistence diagram a vector using the above procedure is continuous when the metric on $\frak{B}$ is the $1$-Wasserstein distance. 
\end{Theorem} 
\begin{Remark} {\em In \cite{images}, estimates proving this result are given both in the case of a general choice of $\phi$ and the special case where $\phi$ is given by Gaussian distributions with fixed variance.  Of course, the estimates in the latter case are stronger .   }
\end{Remark} 
\section{Generalized Persistence} \label{generalized}
\subsection{Introduction}
Persistent homology operates on functors $F$  from the category $\bb{R}_+$ to simplicial complexes, by composing them with the homology functor to obtain a persistence vector space.     It is   useful to consider other parameter categories for the diagrams, which can also help clarify the structure of data sets.  There are at least two such constructions that have been discussed.  The first is {\em zig-zag} persistence, introduced in \cite{zigzag}, and the second is {\em multidimensional persistence}, discussed for example in \cite{multid}.  The first is designed to study the relationship between homology of constructions that are not nested within each other, such as distinct samples from a given space, and the second provides invariants of situations where it is natural to study filtrations of spaces involving more than one parameter, such as filtering by both the scale parameter $R$ and a measure of density.  We describe both extensions of the standard persistent homology methods.  

\subsection{Zig-zag Persistence} Consider the triangulation of $\bb{R}$ whose vertices are the integers. The set of vertices of  the barycentric subdivision of this simplicial complex  is equipped with a partial ordering, by recognizing that its elements are in one to one correspondence with the simplices in the original triangulation, and that that set is equipped with a partial ordering by treating it as a subset of the power set of $\bb{R}$.  In concrete terms, we may view it as in one to one correspondence with the set all integers and half integers, with every integer $n$ being less than or equal to the elements $n \pm \frac{1}{2}$.  We'll refer to this partially ordered set as $\frak{Z}$.  A partially ordered set (and its corresponding category) is said to be connected if any two objects can be connected by a zig-zag sequence of morphisms.  Connected partially ordered subsets of $\frak{Z}$  are always determined by a pair $(x,y)$, where $x$ and $y$ are both integers or half integers, via the rule that assigns to the pair $(x,y)$ the collection of objects $z$ for which $x \leq z \leq y$ (in the total ordering on $\bb{R}$).   

\begin{Definition}
For any category $\underline{C}$, a {\em zig-zag persistence object} in $\underline{C}$ is a functor from  a connected subcategory of $\frak{Z}$ to $\underline{C}$. Suppose further that $\ul{C}$ is equipped with an object $c_0$ that is both initial and terminal.   Then for any connected subcategory $\frak{Z}_0 \subseteq \frak{Z}$ and object $c \in \underline{C}$, we define the {\em interval object} for $\frak{Z}_0 $ and $c$ to be the functor $F = F_{\frak{Z}_0,c}: \frak{Z} \rightarrow \ul{C}$ defined on objects by  $F(x) = c$ for all $x \in \frak{Z}_0$ and $F(x) = c_0$ for any $x \notin \frak{Z}_0$, and on morphisms by $F(x \leq y) = id_{c}$ whenever $x,y \in \frak{Z}_0$.  The behavior on morphisms into or out of $\frak{Z}_0$ is uniquely determined by the fact that $c_0$ is both initial and terminal. When $\ul{C}$ is the category of vector spaces over a field $K$, then it is understood that $c$ will be chosen to be a one-dimensional vector space over $K$. 
\end{Definition}
Given a  zig-zag persistence object in the category of simplicial complexes, we may apply a homology functor $H_i(-;A)$ for an abelian group $A$ to obtain a zig-zag persistence object in the category of abelian groups.  The category of abelian groups has the zero object as an object which is both initial and terminal, and so the notion of interval objects makes sense.  Of course if $A = K$, where $K$ is a field, then we obtain a zig-zag persistence vector space.  It turns out that they can be classified up to isomorphism. 
\begin{Theorem} \label{classify}{\bf (P. Gabriel, \cite{gabriel})}  Let $\frak{Z}\pr\subset \frak{Z}$ be any finite connected subcategory, and let $F$ denote any zig-zag persistence object in the category of finite dimensional  vector spaces over a field $K$  defined on $\frak{Z}\pr$.  Then there is a finite direct sum decomposition 
$$  F \cong \bigoplus _{i} F_i
$$
where each $F_i$ is an interval object for $\frak{Z}_0$ and $K$, where $\frak{Z}_0\subseteq \frak{Z}\pr$ is a connected subcategory of $\frak{Z}\pr$.  Moreover, the sum is unique up to isomorphism and reordering of the sum.  
\end{Theorem} 
\begin{Corollary} The classification of zig-zag persistence vector spaces based on a connected subcategory $\frak{Z}$ is given by barcodes where the intervals have endpoints integers or half integers.  We'll refer to these barcodes as the  zig-zag persistence barcodes of the zig-zag persistence vector space. 
\end{Corollary}
This theory with applications is discussed in \cite{zigzag} and \cite{zzfunction}.
Here are some particular situations in which this construction can be used. 
\begin{enumerate}
\item{{\bf Samples:}}  Given a finite metric space, one can ask to what extent the persistent homology is captured on smaller samples from the data set.  For example, suppose that we have taken a very large uniform sample  $X$ from the circle, and equip them with a metric by restricting  the intrinsic metric (say) on the circle to these points.  We will then with high probability find that the one-dimensional persistence barcode for the Vietoris-Rips construction on $X$ will contain one long bar and many much shorter ones. Supposing that we do not actually know that the sample is coming from a circle, but simply observe that we obtain a one-dimensional barcode with one long bar and many shorter ones.  A hypothesis suggested by this observation is that the data is obtained by sampling from a space with the homotopy type of the circle, but we may wonder if instead it has somehow appeared ``by accident".  One way to provide confirmation of our hypothesis would be to observe that we obtain the same result for various subsamples of our space, and that they are compatible in an appropriate sense.  Zig-zag persistence provides a way for carrying this out.  We suppose that we have chosen samples $X_1, \ldots , X_n \subseteq X$, and create a persistence object in the category of finite metric spaces and distance non-increasing maps 
{$$\divide\dgARROWLENGTH by 4
\begin{diagram}
\node{}  \node{X_1 \cap X_2} \arrow{se}\arrow{sw}\node{}
 \node{\cdot \cdot 
\cdot } \arrow{se} \arrow{sw} \node{}\node{X_n \cap X_{n-1}} \arrow{se} \arrow{sw}\\
\node{X_1} \node{} 
\node{X_2} \node{} 
 \node {X_{n-1}} \node{} \node{X_n}
\end{diagram} 
$$}
If we apply ${\cal V}(-, r) $ for a fixed choice of $r$, guided by the beginning and endpoints of the observed long bar in the barcode for  $X$, we obtain a persistence object in the category of simplicial complexes.  If we further apply $H_1(-,K)$ for a field $K$, we obtain a zig-zag persistence $K$-vector space.  By the classification Theorem \ref{classify} above, we obtain a  decomposition of the resulting zig-zag persistence $K$-vector space.  The  interpretation of the informal idea that each sample has a one-dimensional homology class and that they are consistent is the presence of an interval $K$-vector space for a relatively long interval within the set $\{ 1,\frac{3}{2}, 2, \ldots n - \frac{1}{2}, n \}$, or equivalently a relatively long bar in the zig-zag persistence barcode.  
\item{{\bf Functions on spaces:}} Suppose that we have a topological space $X$  equipped with a continuous map $f: X \rightarrow \bb{R}_+$.  Then we have the ordinary persistence $K$-vector spaces $\{H_i(f^{-1}([0,r],K) \}_r$, which encode information about the evolution of the homology of the sublevel sets of $f$ as $r $ increases.  However, one  might be interested in gaining information about approximations to {\em level sets} instead.  They can provide more useful invariants in a number of cases, and are approachable through zig-zag persistent homology.  We construct a zig-zag diagram of topological spaces as follows.

$$\xy
\xymatrix@C=4pt{ & {f^{-1}(1) } \ar[dl] \ar[dr]&   &{\cdot \cdot \cdot} \ar[dl] \ar[dr]& &f^{-1}(N-1) \ar[dl] \ar[dr]& \\
              f^{-1}I_0& &  f^{-1}I_1 && f^{-1}I_{N-2} & &f^{-1}I_{N-1}}
\endxy
$$
where $I_k = [k,k+1]$  for all $k$.  Again we can apply $H_i(-,K)$ to this diagram, an obtain  a zig-zag persistence barcode.  It contains information about how the spaces $f^{-1}[k,k+1]$ change as $k$ changes, and how they assemble together.  This situation is studied in \cite{zzfunction}.  
\item{{\bf Witness complexes:}} One of the problems with the witness complex is that we have very little theory about the extent to which it reflects accurately the persistent homology of the underlying metric space.  A related problem is that there is no direct relationship between the construction for two different landmark sets.  Even if ${\cal L}_1 \subseteq {\cal L}_2$, there are no maps directly relating $W(X,{\cal L}_1,R)$ and $W(X,{\cal L}_2,R)$ for two different landmark sets.   One approach is to attempt to assess in some manner how consistent the results of the constructions based on ${\cal L}_1$ and ${\cal L}_2$ are.  It turns out that given two landmark sets ${\cal L}_1$ and ${\cal L}_2$, it is possible to construct an intermediate bivariant construction $W(X,\{{\cal L}_1,{\cal L}_2 \},R)$ for which there is an evident diagram 
$$ \begin{diagram}
\node{W(X,{\cal L}_1,R)}  \node{W(X,\{{\cal L}_1,{\cal L}_2 \},R)} \arrow{w} \arrow{e} 
\node{W(X,{\cal L}_2,R)}
\end{diagram} 
$$
This is itself a short zig-zag diagram of length three, but if we have landmark sets $\{ {\cal L}_i \}_{i=0}^N$ we can clearly construct a longer diagram that includes the bivariant constructions $W(X,\{ {\cal L}_i , {\cal L}_{i+1} \}, R)$ for $i=0, \ldots , N-1$.  The construction is quite simple.  Its vertex set is ${\cal L}_1 \times {\cal L}_2$, and a subset 
$$\{ l_1^1 \times l_2 ^1, l_1^2 \times l_2 ^2, \ldots , l_1^k \times l_2^k \} \subseteq {\cal L}_1 \times {\cal L}_2$$ spans a simplex in $W(X,\{ {\cal L}_i , {\cal L}_{i+1} \}, R)$ if and only if there is a point $x \in X$ so that $x$ is a witness  for $\{ l_1^1, l_1^2, \ldots , l_1^k \}$ and $\{l_2^1, l_2 ^2, \ldots , l_2^k \}$ in the complexes $W(X,{\cal L}_1 ,R)$ and  $W(X,{\cal L}_1 ,R)$, respectively. The projections to  $W(X,{\cal L}_1 ,R)$ and  $W(X,{\cal L}_1 ,R)$ are defined in the evident way.  

\end{enumerate} 
\subsection{Multidimensional Persistence}
There are many situations where it can be useful to introduce families of spaces varying with more than one real parameter.  For example, in any kind of topological analysis of data sets, it usually is the case that if one considers the persistent homology of the entire data set, the presence of outliers means that we do not typically obtain the``right homology".  For example, if we have data sampled from the unit  circle, but a small number of points sprinkled throughout the unit disc, then persistent homology will end up reflecting the homology of the disc rather than that of the circle.  This is often circumvented by selecting only the points of sufficient density, as measured by a density estimator, since outliers will typically have very low density.  The question then becomes, though, how to choose the threshold for density.  Also, it turns out that in general, there will be variation through different topologies as one changes the threshold.  The solution to this problem is to attempt to study all thresholds at once, just as we do when considering the scale parameter in the Vietoris-Rips construction.   This leads us to the following definition.
\begin{Definition}  Let $\ul{C}$ be a category. Then by an $n$-dimensional persistence object in $\ul{C}$ we mean a functor $\bb{R}_+^n \rightarrow \ul{C}$, where $\bb{R}_+^n$ denotes the $n$-fold product of copies of the category $\bb{R}_+$.   
\end{Definition} 
\begin{Example} {\em Let $X$ be any metric space, and suppose that $X$ is equipped with a function 
$f:X \rightarrow \bb{R}$.  It might be a density estimator, but it might also be a measure of centrality.  Then if $X[s]$ denotes the subset $\{ x \in X | f(x) \leq s \}$, we obtain family of spaces ${\cal V}(X[s],r)$, and by applying homology with coefficients in a field $K$, we obtain a $2$-dimensional 
persistence $K$-vector space parametrized by the pair $(r,s)$. While density is used as described above to remove outliers or noise, the case of a centrality measure allows one to capture the presence of the analogues of ends in a finite metric space.  }  \end{Example} 
\begin{Example}{\em  Multidimensional persistent homology can also be used to capture geometric information that is not strictly topological (see \cite{shapes}). For example, given a Riemannian manifold, one can compute Gaussian curvature at each point and filter by that quantity.  Considering the entire manifold, one can   obtain a one-dimensional persistence vector space by applying homology over a field $K$.  For computational purposes, though, we would need to deal with a sample and use a second parameter, namely the scale parameter in a Vietoris-Rips complex.  This kind of analysis can for example be used to distinguish between various ellipsoids.  
}
\end{Example} 
The equivalence of categories in Proposition \ref{pclassify} has the following  straightforward analogue.  
\begin{Proposition} \label{multiequivalence} The category of $n$-dimensional persistence vector spaces over $K$ is equivalent to the category of $\bb{R}_+^n$-graded modules over the $\bb{R}^n_+$-graded ring $K[\bb{R}_+^n]$. The analogous result with $\bb{R}_+^n$ replaced by $\bb{N}^n$ also holds.   
\end{Proposition}

Although useful, this result does not give us a classification of multidimensional persistence $K$-vector spaces analogous to the barcode classification that works in the one-dimensional case.  The reason can be understood as analogous to the commutative algebraic situation, where finitely generated  $K[x]$-modules can be classified  because $K[x]$ is a principal ideal domain, but  $K[x_1, \ldots , x_n]$-modules cannot be classified for $n \geq 2$.   In a sense it is provable that there is no classification strictly analogous to the $n = 1$ case, because the classification in the case $n \geq 2$ depends on the structure of the field $K$, as is demonstrated in \cite{multid}.  This is not the case when $n = 1$, because the classification is always by barcodes, and that is independent of the structure of $K$.  

Although there is no complete classification  of multidimensional persistence vector spaces, there exist interesting invariants.  For any $k$-dimensional  persistence $K$-vector space $\{V_x \}_{x \in \bb{R}^k}$, and  pair of points $x,y \in \bb{R}^k$, with $x \leq y$ in the natural partial ordering on $\bb{R}^k$, we can define $r(x,y)$ to be the rank of the linear transformation $ V_x \rightarrow V_y$.  We extend the definition to all pairs $x,y \in \bb{R}^k$ by setting $r(x,y) = 0$ when $x$ is not less than or equal to $y$.  The function $r$ is therefore defined on $\bb{R}^k \times \bb{R}^k$, and we refer to it as the {\em rank invariant}.  In the case $k=1$, the rank invariant is complete, in that it differentiates between distinct barcodes. There is an analogue for multidimensional persistence $K$-vector spaces of  the results described in   Section \ref{algebraic}.  To state it, we first define an analogue of barcodes.  By a {\em cube} in $\bb{R}^k$, we mean a set of the form $I_1 \times I_2 \times \cdots \times I_k$, where each $I_s$ is an interval $[a_s,b_s]$, and write $C(a_1, \ldots , a_k, b_1, \ldots , b_k)$ for this cube.  We write ${\cal C}_k$ for the set of all $k$-dimensional cubes.   There is a straightforward analogue to the space $\frak{B} $ defined in Section \ref{featureintro}.  We first define $B_n = Sp^n({\cal C}_k)$, and define $\frak{B}(k)_+= \coprod Sp^n({\cal C}_k)$.  Next, we say a cube $C(a_1, \ldots , a_k, b_1, \ldots , b_k)$ is {\em negligible} if $a_i = b_i $ for some $i$, and define an equivalence relation $\sim$ on $\frak{B}(k)_+$ to be the equivalence relation generated by the relation 
$$ \{ C_1, C_2, \ldots C_n \} \sim \{ C_1, \ldots , C_{n-1} \} \mbox{ for all negligible } C_n 
$$  
We define $\frak{B}(k)$ to be $\frak{B}(k) _+/\sim$.  Every equivalence class under $\sim$ has a unique minimal representative consisting entirely of non-negligible cubes.  Let ${\cal M}(k)$ denote the set of isomorphism classes of $k$-dimensional persistence $K$-vector spaces.  For each cube $C = C(a_1, \ldots , a_k, b_1, \ldots , b_k )$ we let $\mu (C)$ denote the isomorphism class of the the $k$-dimensional persistence $K$-vector space $\{ V_{\vec{x}}\}_{\vec{x} \in \bb{R}^k_+}$ for which $V_{\vec{x}} = K$ whenever $\vec{x} \in C$, $V_{\vec{x}} = \{ 0 \}$ whenever $\vec{x} \notin C$, and for which all induced morphisms $V_{\vec{x}} \rightarrow V_{\vec{y}}$ for $\vec{x} \leq \vec{y}$ and $\vec{x},\vec{y} \in C$ are equal to the identity.  There is an obvious map $\theta: \frak{B}(k) \rightarrow {\cal M}(k)$ which assigns to each minimal representative $\{C_1, \ldots , C_n \}$ the direct sum $\oplus_i \mu(C_i)$.  

\begin{Theorem} The constructions above satisfy the following properties.  
\begin{enumerate}
\item{The set $\frak{B}(k)$ is a subset of the set  of real points of an affine scheme $Spec(A)$. }
\item{The ring $A$ is complicated, but there is a  $\Bbb{G}_m$-action on $Spec(A)$ that allows us to define  a more manageable subring $A^{fin} \subseteq A$.  }
\item{$A^{fin}$ is isomorphic to the polynomial ring $\bb{R}[x_{\vec{a},\vec{b}}]$ where $\vec{a}$ and $\vec{b}$ are $k$-vectors of integers for which $a_i \geq 1$ and $b_i \geq 0$ for all $i$. }
\item{The ring $A^{fin}$ separates points in $\frak{B}(k)$, and maps injectively to the ring of all real-valued functions on $\frak{B}(k)$.  }
\item{There is a natural lift  of the ring homomorphism $A^{fin} \rightarrow F(\frak{B}(k),\bb{R})$ along $\theta$  to a ring homomorphism $j: A^{fin} \rightarrow F({\cal M}(k), \bb{R})$. $F(X,\Bbb{R})$ denotes the ring of real valued functions on a set $X$.    }
\item{ For any $\alpha \in A^{fin}$, the function $j(\alpha): {\cal M}(k) \rightarrow \bb{R}$ factors through the rank invariant.  For any two elements $X,Y \in {\cal M}(k)$,  if  for all $\alpha \in A^{fin}$ , $j(\alpha )(X) = j(\alpha) (Y)$,  then  the rank invariants of $X$ and $Y$ are equal.   }
\end{enumerate} 
\end{Theorem} 

This result gives one approach to the study of multidimensional persistence.  There is a great deal of other work on this topic.  The paper \cite{chacholski2} deals with computing persistent homology using commutative algebra techniques, via the equivalence of categories from Proposition \ref{multiequivalence}.  They demonstrate that the multigrading yields significant simplification.  In \cite{chacholski1} and \cite{chacholski2}, an algebraic  framework is constructed for dealing with the fact that there is alway noise in the applications to data analysis.  In \cite{ferri}, multidimensional persistence is studied by examining the family of all one dimensional persistence modules obtained by considering lines with varying angles in the persistence domain.  M. Lesnick in \cite{lesnickthesis} has defined metric properties of the set of isomorphism classes of multidimensional persistence vector spaces, and proven uniqueness results for them. Finally, in \cite{lesnickwright}, software is developed for visualization and interrogation of two-dimensional persistence vector spaces.  

One generalization that has not been studied yet is to multidimensional persistence where some of the persistence directions might be ``zig-zag" directions rather than ordinary persistence directions.  Formally, this would mean functors from the categories of the form $\bb{R}_+ ^m \times \frak{Z}^n$.  This would be very useful in a number of situations.  For example, in the zig-zag constructions discussed in Section \ref{generalized} for samples and for witness complexes, we were forced to choose a threshold for the scale parameter.  If we had a way of representing functors from $\bb{R}_+ \times \frak{Z}$ to vector spaces, we would not be forced to make this selection.  

\section{Coverage and Evasion Problems}
There is an interesting set of technologies used for sensing of various kinds called {\em sensor nets}.  A sensor net consists of a collection of  sensors distributed throughout a domain. The sensors are very primitive in the sense that they are only capable of sensing the presence of an intruder or of another sensor within a fixed detection radius $R$. We also assume each sensor is given an identifying label or number, and that other sensors can sense that identifier when they are within the radius $R$ of each other.         One problem of interest is whether or not the balls of radius $R$ cover the domain, and it does not have an immediate solution due to the fact that the positions of the sensors are not available.   V. De Silva, and R. Ghrist have developed a very interesting method for addressing this problem based on persistent homology (see \cite{desilva1} and \cite{desilva2}).  The rough idea is as follows.  Suppose that one has a domain $D$ in the plane, with a  connected and compact curve boundary $\partial D$, that  one has a collection of points $\{v_i \}_{i \in I}$ in the region, one for each sensor, and that one knows in some way that $\partial D$ is covered by the open balls   $B_R(v_i)$. Let ${\cal U}$ denote the family of open subsets  $\{B_R(v_i) \cap D\}_{i \in I}$ of $D$, and let ${\cal U}^{\partial}$ denote the covering   $\{B_R(v_i) \cap \partial D\}_{i \in I}$ of $\partial D$.  Suppose further that one knows that spaces 
$$ B_R(v_{i_1}) \cap \cdots \cap B_R(v_{i_s}) \cap D \mbox{ and }   B_R(v_{i_1} )\cap \cdots \cap B_R(v_{i_s}) \cap \partial D
$$
are all either contractible or empty, for all choices of subsets $\{i_1, \ldots , i_s \} \subseteq I$.  The conditions assure that $D$ is homotopy equivalent to the nerve of ${\cal U}$, and that $\partial D$ is homotopy equivalent to the nerve of the covering ${\cal U}^{\partial}$, as a consequence of the nerve theorem.    It further assures that the pair $(N_.{\cal U}, N_.{\cal U}^{\partial})$ is equivalent to the pair $(D, \partial D)$.  For any field, the relative group $H_2(D, \partial D;K) \cong K$, since $D$  is a connected orientable manifold with boundary $\partial D$, and it follows that the relative group $H_2(N_.{\cal U}, N_. {\cal U}^{\partial};K) \cong K$.    On the other hand, suppose that the sets $B_R(v_i)$ do not cover all of $D$, and let $D_0 \subseteq D$ denote the union 
$$ \bigcup _{i} B_r(v_i) \cap D
$$
The space $D_0$ is a non-compact  manifold with boundary, and consequently the relative group $H_2(D_0, \partial D;K) \cong 0$.  As above, it will follows that $$H_2(N_.{\cal U}, N_.{\cal U}^{\partial};K) \cong 0$$
Consequently, the simplicial complex of the nerve of the covering, which can be computed using the information available from the sensors, determines whether or not we have a covering based on its simplicial homology.  The conditions on the coverings given above are of course impossible to verify, but DeSilva and Ghrist are able to formulate a persistent homology condition that is a reasonable substitute, and which gives a homological criterion in terms of  a $2$-dimensional persistent homology group which is sufficient to guarantee coverage.  

In order to understand the result in \cite{desilva2}, we first observe that the information from the sensors do not give us access to the \v{C}ech complex, since we have no way of determining the intersection of balls without precise knowledge of the distances between their centers. However, we do have access to the Vietoris-Rips complex, since for any pair of points, we can tell whether or not they are within a distance $R$, where $R$ is the detection radius.  We also have the comparison results for the Vietoris-Rips complex and the \v{C}ech complex given in Proposition \ref{interleave}.  The overall idea in \cite{desilva1} and \cite{desilva2} is to leverage the relationship between the complex we have access to (Vietoris-Rips) and the complex from which we can deduce coverage (\v{C}ech).  In order to formulate such a result, we assume that we are attaching a second number to each sensor, namely its {\em covering radius} $R_c$.  It is understood that each sensor covers a disc of radius $R_c$ around it, and that it can detect other sensors at the detection radius $R$ given above.  We further assume that $R_c \geq R/\sqrt{3}$.  This allows us to guarantee that if $\{ x_0,x_1,x_2 \}$ forms a two simplex in the Vietoris-Rips complex ${\cal V}(X,R)$, then they span a two simplex in $\check{C}(X,R_c)$, by   the second statement in Proposition \ref{interleave}.    There are now the following assumptions made in \cite{desilva2}.  
\begin{enumerate}
\item{The sensors lie in a compact connected domain $D \subseteq \bb{R}^2$, whose boundary $\partial D$ is connected and piecewise linear with vertices called  {\em fence nodes}. }
\item{Each fence node $v$ is within $R$ of its fence node neighbors on $\partial D$}
\end{enumerate}
The following theorem is proved in  \cite{desilva1} and \cite{desilva2}.  
\begin{Theorem} Let ${\cal R}$ denote the Vietoris-Rips complex of the set of all sensors, and let ${\cal F}$ denote the subcomplex on the fence vertices.  If the sensors satisfy the conditions above, and if there exists $[\alpha] \in H_2({\cal R},{\cal F})$ so that $\partial ([\alpha]) \neq 0 $, where 
$\partial: H_2({\cal R},{\cal F}) \rightarrow H_1({\cal F})$ is the connecting homomorphism, then the balls of radius $R_c$ around the sensors cover $U$.  

\end{Theorem} 
This theorem is in some situations not ideal, due to the strong assumptions on the boundary.  In \cite{desilva2}, it is shown that the use of the persistent homology of the pair $({\cal R}, {\cal F})$ can be used to obtain coverage results with much weaker hypotheses.  

Another interesting direction is the study of time varying situations, where the sensors move in time.  In this case, there are situations where the balls around the sensors do not cover the region at any fixed time, but that no ``evader" can avoid being sensed at some time.  This kind of problem is referred to as an {\em evasion problem}, and has been studied in \cite{evasionadams} and \cite{ghrist}.  The two approaches are quite distinct, the approach in \cite{evasionadams} using zig-zag persistence, and the approach in \cite{ghrist} develops a new kind of cohomology with semigroup coefficients.  The  approach in \cite{ghrist} yields  ``if and only if" results.  
\section{Probabilistic Analysis}
\subsection{Random Complexes} A very interesting direction of research is the study of the distributions on the space of barcodes $\frak{B}$ (defined in Section \ref{wholemetrics}) obtained by sampling points from Euclidean space using various models of randomness, i.e. sampling from various distributions on $\bb{R}^n$, or using the theory of random graphs \cite{bollobas}.  Since we do not have a library of well understood distributions on $\frak{B}$, one can instead study the distributions on the real line obtained by pushing forward a distribution on $\frak{B}$ along a map from $\frak{B} $ to $\bb{R}$.     An interesting such map is $$  \{(x_1,y_1), (x_2,y_2), \ldots ,( x_n, y_n) \}  = B  \rightarrow \mbox{max}\{\frac{y_1}{x_1}, \ldots , \frac{y_n}{x_n} \} = \lambda(B)
$$
In fact, this map is only defined for barcodes in $\frak{B}$ for which the endpoints of the intervals in the barcode all lie in $\{ x \in \bb{R}| x >0 \}$.  When one is computing homology in dimensions $\geq 1$  of a \v{C}ech or Vietoris-Rips complex of a set of points in Euclidean space, the barcode satisfies this property.  
A very interesting result in this direction is proved in \cite{kahleetal}.  They proceed by sampling from $[0,1]^d$  using a uniform Poisson process of intensity $n$.  This means that the sampling is done from a uniform distribution on $[0,1]^d$, but that the number of points sampled is a governed by a Poisson distribution.  The description of all these notions is beyond the scope of the present paper, but we refer the reader to \cite{stochprocess}.  The main result of \cite{kahleetal} is the following. 
\begin{Theorem} We suppose that $d$ and $n$ are as above, that $d\geq 2$,   that we are computing the $k$-dimensional persistence barcode $B$, and that $ 1 \leq k \leq d-1$.  Let 
$$ \Delta _k(n) = { \Big (}\frac{\mbox{\em log }n}
{\mbox{{\em log log} n}}{\Big )}^{\frac{1}{k}}$$
Then there exist constants $A_k$ and $B_k$ so that 
$$ \column{\mbox{\em lim}}{n \rightarrow \infty} \bb{P}{\Big (} A_k \leq \frac{\Pi _k(n)}{\Delta _k(n)}\leq B_k {\Big )} = 1
$$
where $\bb{P}$ denotes probability, and where $\Pi _k(n)$ denotes the value of $\lambda (B)$ for a barcode generated as above.    
\end{Theorem}
\subsection{Robust Estimators} The stability theorem in Section \ref{stable}  deals with the effect on persistence barcodes of small perturbations in the metric space, where perturbations are small in the sense of the Gromov-Hausdorff distance.  In reality, though, one expects that in a perturbation of a metric space, a small number of distances may  undergo relatively large perturbations.  However, one believes that the number will be small, and that the points involved will be of small measure in an underlying measure.  In order to deal with this problem, one incorporates a measure-theoretic component in one's definitions.  
\begin{Definition}  By a {\em metric measure space} we will mean a complete separable metric space $M$ equipped with a Borel measure $\mu$.
\end{Definition} 
In \cite{greven}, a metric $d_{GPr}$  called the {\em Gromov-Prokhorov} metric is introduced on the measure preserving  isometry classes of compact metric measure spaces.  It is constructed by combining the Gromov-Hausdorff metric on the isometry classes of compact metric spaces with the Prokhorov metric $d_{Pr}$ on measures on a fixed metric space (see \cite{prokhorov}) by methods which we will not discuss here.  

The paper \cite{blumberg} studies the distributions on the space of persistence barcodes arising from the persistence barcodes obtained by sampling from a fixed metric measure spaces. More precisely, they study distributions on the completion of the metric space $\overline{\frak{B}}$ of $\frak{B}$ equipped with the bottleneck distance.   Let $\mu(n,k,X)$ denote the distribution on $\overline{\frak{B}}$ which arises from sampling a set  $S$ of $n$  points on a metric measure space $X$, and computing the $k$-dimensional barcode on $S$.  We have the following. 
\begin{Theorem} The inequality 
$$ d_{Pr}(\mu(n,k,X),\mu(n,k,X^{\pr})) \leq n d_{GPr}(X,X^{\pr})
$$
holds for all compact metric measure spaces $X$ and $X{\pr}$.  
\end{Theorem} 
This result is the used in \cite{blumberg} to develop {\em robust statistics} for distinguishing the results of sampling from a fixed metric space.  Robust statistics are computable quantities attached to samples from distributions that are relatively insensitive  to small changes in parameter values in the distribution from which the samples are gathered, and also to the presence of outliers.  An elementary example of this idea is the median, which is relatively insensitive to outliers and is considered a robust statistic, while the mean is not.  For the problem at hand, \cite{blumberg} defines a precise notion of robustness. 

\begin{Definition} Let $f$ be a function from the set of isomorphism classes of finite metric spaces to a metric space $(W,d_W)$.  We say the $f$ is {\em robust with robustness coefficient $r > 0$} if for any nonempty finite metric space $(X,d_X)$,  if for any nonempty finite metric space $(X,d_X)$, there exists a bound $\delta$ such that for all isometric embeddings of $X$ in a finite metric space $(X^{\pr}, d_{X^{\pr}})$ for which $\#(X^{\pr})/\#(X) < 1+r$, it is the case that $d(f(X),f(X^{\pr})) < \delta$.  There is a corresponding uniform notion that states that the bound $\delta$ may be chosen universally, for all $X$.
  There is a corresponding uniform notion that states that the bound $\delta$ may be chosen universally, for all $X$.
\end{Definition} 
We now obtain the following result for finite metric spaces, which are being regarded as metric measure spaces by assigning to each metric space the uniform measure. 

\begin{Theorem}{\em (\cite{blumberg})} For fixed $n$ and $k$, $\mu(n,k,-)$ is uniformly robust with robustness coefficient $r$ and estimate bound $\delta=nr/(1+r)$ for any $r$.  
\end{Theorem} 

Since the space $\overline{\frak{B}}$ is relatively inaccessible, it is useful to use this result to construct real valued statistics that also satisfy the robustness property.  One way to do this is to consider a fixed reference distribution ${\cal P}$ on $\overline{\frak{B}}$, and define  $\Delta _{\cal P}(n,k,X)$ to be the number $d_{Pr}(\mu(n,k,X), {\cal P})$.  

\begin{Theorem}{\em (\cite{blumberg})} For fixed $n$ and $k$, and ${\cal P}$, $\Delta _{\cal P}(n,k,-)$ is uniformly robust with robustness coefficient $r$ and estimate bound $nr/(1+r)$ for any $r$.  
\end{Theorem} 
It is also possible to obtain a somewhat simpler result, which does not require calculation of the full distribution $\mu(n,k,X)$.  Instead of fixing a reference distribution ${\cal P}$, we choose a reference barcode $B \in \overline{\frak{B}}$, and  define $\Delta _{B}^{med}(n,k,X) $ to the median of the distribution of $d_B(B,-)$ applied to samples  of $k$-dimensional barcodes attached to samples of size $n$ taken from the metric measure space $X$.  
\begin{Theorem} For fixed $n$, $k$, and $B$, the function $\Delta _{B}^{med}(n,k,-)$ from finite metric spaces (with uniform probability measure) is robust with robustness coefficient greater than $\mbox{\em ln }2/n$.  
\end{Theorem}

\subsection{Random Fields}
Suppose that we have a function $f$ on a manifold.  We have seen in Section \ref{wholemetrics} that we can associate to $f$ the persistence vector spaces $\{ H_i(f^{-1}((- \infty, r],K) \}_r$ for $i$ a non-negative integer and $K$ a field, and further that there are stability results that show that small changes in the function $f$ lead to small changes, as measured by the bottleneck distance, in the corresponding sublevel set barcode.  One can also ask, though, what the expected behavior of various of the features attached to barcode space is for a class of functions chosen at random.  An initial question is what one means  by a function chosen at random.  The notion of a {\em random field} is defined in \cite{adler} as follows.  
\begin{Definition}  By a {\em real-valued random field} on a topological space $T$ we mean a measurable mapping 
$$F: \Omega \rightarrow \bb{R}^T
$$ 
where $\bb{R}^T$ denotes the set of all real-valued functions on $T$ and $(\Omega, {\cal F}, \bb{P})$ is a complete probability space.   $F$ creates a probability measure on $\bb{R}^T$, from which one can sample.  Similarly, one can define an {\em $s$-dimensional vector-valued  random field} as a measurable mapping 
$$F: \Omega \rightarrow (\bb{R}^T)^s
$$ 
\end{Definition} 

 The idea here is that rather than being a function, a random field is an assignment to each $t \in T$ a distribution on $\bb{R}$, rather than a fixed value.   Each of the restrictions $\bb{R}^T \rightarrow \bb{R}^{\{ t \}} \cong \bb{R}$ produces a random variable, and therefore the corresponding distribution, which we denote $F_t$.  In fact, for any finite set of points $t_1, \ldots , t_n \in T$, we obtain a distribution on $\bb{R}^n$ which we denote $F_{t_1, \ldots , t_n }$.   There is a particular class of random fields called the {\em Gaussian random fields} that is particularly amenable to analysis.  
 
 \begin{Definition} A real valued random field $F$ is a {\em Gaussian random field} if for all $n$,  the $n$-dimensional distributions $F_{t_1, \ldots , t_n }$ are multivariate Gaussian distributions on $\bb{R}^n$.  An $s$-dimensional  random field is Gaussian, if all the distributions $F_{t_1,\ldots , t_n }$ are multivariate Gaussian distributions on $\bb{R}^{n \times s}$.  Note that in Gaussian fields, the behavior of the random function is completely determined by the expectation function $m(t) = E(F_{t_1, \ldots , t_n })$, and the covariance function $C$.    
 \end{Definition} 
 
The first example of this construction comes out of work of  Wiener (see \cite{wiener} and  \cite{baldi}), using analysis of Brownian motion.  Wiener studied the case $T = \bb{R}_+$, and produced a Gaussian random field $W$, where the expected value of $W_t$ is always $=0$ and the variance is given by $C(s,t) = \mbox{min}(s,t)$.    He also showed that when one samples from the associated measure on $\bb{R}^{\bb{R}_+}$, one obtains continuous functions with probability 1, and so one calls $W$ a {\em continuous Gaussian field}. 
Given any analytic property of functions on manifolds, such as $k$-th order differentiability, smoothness, or the property of being a Morse function, one can create and study Gaussian random fields whose samples have the given property with probability 1.  Further, there are frequently a priori conditions on the covariance function of the random field that can be readily verified, and guarantee the satisfaction of such properties.  

The paper \cite{bobrowski} proves a result concerning the persistent homology of the sublevel sets of functions sampled from Gaussian random fields.  We consider the real valued function $\sigma$  on barcodes given by 
$$ \sigma  \{ [a_1,b_1], \ldots , [a_n,b_n] \} = \sum _i (b_i - a_i)
$$
For any fixed $x \in \bb{R}$ and barcode $\beta = \{ [a_1,b_1], \ldots , [a_n,b_n] \}$, we define the {\em $x$-truncation of $\beta$},  $\beta [x]$,  to be the barcode 
$$\{ [a_1, \mbox{min}(b_1,x)], \ldots , [a_n, \mbox{min}(b_n,x)] \}
$$
where it is understood that for any $i$ such that $x \leq a_i$, the interval $[a_i,b_i]$ is simply deleted.  
Finally, we define  
$$\chi^{pers}(M,f, x) = \sum _{i=0}^{\infty} (-1)^i \sigma (\{ H_i(f^{-1}((- \infty, r],K) \}_r[x])
$$
In \cite{bobrowski}, the following result is proved concerning the distribution of $\chi ^{pers}(M,f,x)$ for Gaussian random fields on Riemannian manifolds which produce Morse functions with probability one. 

\begin{Theorem} Let $M$ be a closed $d$-dimensional Riemannian manifold, and let $F$ be a smooth real valued Gaussian random field so that $F$ is Morse with probability 1, and so that the mean is identically zero and the variance is identically equal to $1$.  Then for any $x \in \bb{R}$, we have 
$$\bb{E}\{ \chi^{pers}(M,f,x) \} = \chi(M)(\varphi(x) + x \Phi (x)) $$ $$+ \varphi(x)  \sum _{j=1}^d (2 \pi ) ^{-j/2} {\cal L} _j (M)H_{j-2}(-x)
$$
where
\begin{enumerate}
\item{$H_n$ denotes the $n$-th Hermite polynomial.}
\item{$\varphi(x) = (2 \pi )^{-1/2}e^{-x^2/2}$ is the density function for the standard Gaussian distribution.}
\item{$\Phi (x) = \int _{- \infty}^x \varphi (u) du$}
\item{${\cal L}_j(M)$ denotes the $j$-th Lipschitz-Killing curvature of $M$ (defined for example in \cite{adler}, Section 7.6) with respect to a metric constructed from the covariance metric attached to $F$.   }
\end{enumerate} 
\end{Theorem} 
\begin{Remark} {\em The point of this result is that it gives a theoretical estimate for the persistent Euler characteristic of sublevel sets in terms of classical invariants of the manifols.  Also, the result in \cite{bobrowski} is actually proved in a much more general context, that of regular stratified spaces and stratified Morse theory, which in particular permits the study of manifolds with boundary. It also includes the study of random fields that are of the form $G \compcirc (F_1, \ldots , F_k)$, where $G$ is a deterministic function from $\bb{R}^k $ to $\bb{R}$, and $(F_1, \ldots , F_k )$ is a vector-valued Gaussian random field.  }
\end{Remark}


\end{document}